% To create a table of content, Tex this. Then move toc.tmp to toc.tex.
% Then uncomment "\listtoc" below.

%%%%%%%%%%%%%%%%%%  tex macros for preprints, cm version %%%%%%%%%%%%%%
%                     (P. Ginsparg, last updated 9/91)
%                if confused, type `b' in response to query 
%
%---------------------------------------------------------------------%
%% site dependent options: 
%% \unredoffs and \redoffs define horizontal and vertical offsets 
%% respectively for unreduced and reduced modes. \speclscape defines
%% the \special{} call that sets printer to landscape (sideways) mode.
%% from standard set below, leave uncommented as appropriate or redefine
%
%%% next 400dpi
%\def\unredoffs{} \def\redoffs{\voffset=-.31truein\hoffset=-.48truein}
%\def\speclscape{\special{landscape}}
%
%%% apple lw
\def\unredoffs{} \def\redoffs{\voffset=-.31truein\hoffset=-.59truein}
\def\speclscape{\special{ps: landscape}}
%
%%% qms lasergrafix:
%\def\unredoffs{} \def\redoffs{\voffset=-.4truein\hoffset=.125truein}
%\def\speclscape{\special{qms: landscape}}
%
%%% saclay A4 paper:
%\def\unredoffs{\hoffset-.14truein\voffset-.2truein} 
%\def\redoffs{\voffset=-.45truein\hoffset=-.21truein} 
%\def\speclscape{\special{landscape}}
%
%---------------------------------------------------------------------%
%
\newbox\leftpage \newdimen\fullhsize \newdimen\hstitle \newdimen\hsbody
\tolerance=1000\hfuzz=2pt
\catcode`\@=11 % This allows us to modify PLAIN macros.
%\def\bigans{b }
%\message{ big or little (b/l)? }\read-1 to\answ
%
\ifx\answ\bigans\message{(This will come out unreduced.} %\let\l@r=L
\magnification=1200\unredoffs\baselineskip=16pt plus 2pt minus 1pt
\hsbody=\hsize \hstitle=\hsize %take default values for unreduced format
\else\message{(This will be reduced.} \let\l@r=L
\magnification=1000\baselineskip=16pt plus 2pt minus 1pt \vsize=7truein
\redoffs \hstitle=8truein\hsbody=4.75truein\fullhsize=10truein\hsize=\hsbody
\output={\ifnum\pageno=0 %%% This is the HUTP version
	\shipout\vbox{\speclscape{\hsize\fullhsize\makeheadline}
	  \hbox to \fullhsize{\hfill\pagebody\hfill}}\advancepageno
	\else
	\almostshipout{\leftline{\vbox{\pagebody\makefootline}}}\advancepageno 
	\fi}
\def\almostshipout#1{\if L\l@r \count1=1 \message{[\the\count0.\the\count1]}
      \global\setbox\leftpage=#1 \global\let\l@r=R
 \else \count1=2
  \shipout\vbox{\speclscape{\hsize\fullhsize\makeheadline}
      \hbox to\fullhsize{\box\leftpage\hfil#1}}  \global\let\l@r=L\fi}
\fi
%---------------------------------------------------------------------
%
\newcount\yearltd\yearltd=\year\advance\yearltd by -1900

\def\Title#1#2{\nopagenumbers\abstractfont\hsize=\hstitle\rightline{#1}%
\vskip 1in\centerline{\titlefont #2}\abstractfont\vskip .5in\pageno=0}
\def\Date#1{\vfill\leftline{#1}\tenpoint\supereject\global\hsize=\hsbody%
\footline={\hss\tenrm\folio\hss}}% 	restores pagenumbers
%
%       use following instead of \Date on the preliminary draft, 
%       puts date/time on each page in big mode, writes labels in margins

\def\draftmode{\message{ DRAFTMODE }\def\draftdate{{\rm preliminary draft:
\number\month/\number\day/\number\yearltd\ \ \hourmin}}%
\headline={\hfil\draftdate}\writelabels\baselineskip=20pt plus 2pt minus 2pt
 {\count255=\time\divide\count255 by 60 \xdef\hourmin{\number\count255}
  \multiply\count255 by-60\advance\count255 by\time
  \xdef\hourmin{\hourmin:\ifnum\count255<10 0\fi\the\count255}}}
%       use \nolabels to get rid of eqn, ref, and fig labels in draft mode
\def\nolabels{\def\wrlabeL##1{}\def\eqlabeL##1{}\def\reflabeL##1{}}
\def\writelabels{\def\wrlabeL##1{\leavevmode\vadjust{\rlap{\smash%
{\line{{\escapechar=` \hfill\rlap{\sevenrm\hskip.03in\string##1}}}}}}}%
\def\eqlabeL##1{{\escapechar-1\rlap{\sevenrm\hskip.05in\string##1}}}%
\def\reflabeL##1{\noexpand\llap{\noexpand\sevenrm\string\string\string##1}}}
\nolabels
%
% tagged sec numbers
\global\newcount\secno \global\secno=0
\global\newcount\meqno \global\meqno=1
\def\newsec#1{\global\advance\secno by1\message{(\the\secno. #1)}
%\ifx\answ\bigans \vfill\eject \else \bigbreak\bigskip \fi  %if desired
\global\subsecno=0\eqnres@t\noindent{\bf\the\secno. #1}
\writetoca{{\secsym} {#1}}\par\nobreak\medskip\nobreak}
\def\eqnres@t{\xdef\secsym{\the\secno.}\global\meqno=1\bigbreak\bigskip}
\def\sequentialequations{\def\eqnres@t{\bigbreak}}\xdef\secsym{}
\global\newcount\subsecno \global\subsecno=0
\def\subsec#1{\global\advance\subsecno by1\message{(\secsym\the\subsecno. #1)}
\ifnum\lastpenalty>9000\else\bigbreak\fi
\noindent{\it\secsym\the\subsecno. #1}\writetoca{\string\quad 
{\secsym\the\subsecno.} {#1}}\par\nobreak\medskip\nobreak}
\def\appendix#1#2{\global\meqno=1\global\subsecno=0\xdef\secsym{\hbox{#1.}}
\bigbreak\bigskip\noindent{\bf Appendix #1. #2}\message{(#1. #2)}
\writetoca{Appendix {#1.} {#2}}\par\nobreak\medskip\nobreak}
%
%       \eqn\label{a+b=c}	gives displayed equation, numbered
%				consecutively within sections.
%     \eqnn and \eqna define labels in advance (of eqalign?)
%
\def\eqnn#1{\xdef #1{(\secsym\the\meqno)}\writedef{#1\leftbracket#1}%
\global\advance\meqno by1\wrlabeL#1}
\def\eqna#1{\xdef #1##1{\hbox{$(\secsym\the\meqno##1)$}}
\writedef{#1\numbersign1\leftbracket#1{\numbersign1}}%
\global\advance\meqno by1\wrlabeL{#1$\{\}$}}
\def\eqn#1#2{\xdef #1{(\secsym\the\meqno)}\writedef{#1\leftbracket#1}%
\global\advance\meqno by1$$#2\eqno#1\eqlabeL#1$$}
%
%			 footnotes
\newskip\footskip\footskip14pt plus 1pt minus 1pt %sets footnote baselineskip
\def\footnotefont{\ninepoint}\def\f@t#1{\footnotefont #1\@foot}
\def\f@@t{\baselineskip\footskip\bgroup\footnotefont\aftergroup\@foot\let\next}
\setbox\strutbox=\hbox{\vrule height9.5pt depth4.5pt width0pt}
\global\newcount\ftno \global\ftno=0
\def\foot{\global\advance\ftno by1\footnote{$^{\the\ftno}$}}
%
%say \footend to put footnotes at end
%will cause problems if \ref used inside \foot, instead use \nref before
\newwrite\ftfile   
\def\footend{\def\foot{\global\advance\ftno by1\chardef\wfile=\ftfile
$^{\the\ftno}$\ifnum\ftno=1\immediate\openout\ftfile=foots.tmp\fi%
\immediate\write\ftfile{\noexpand\smallskip%
\noexpand\item{f\the\ftno:\ }\pctsign}\findarg}%
\def\footatend{\vfill\eject\immediate\closeout\ftfile{\parindent=20pt
\centerline{\bf Footnotes}\nobreak\bigskip\input foots.tmp }}}
\def\footatend{}
%
%     \nref\label{text}
% generates a number, assigns it to \label, generates an entry.
% To list the refs on a separate page,  \listrefs
%
\global\newcount\refno \global\refno=1
\newwrite\rfile
\def\ref{[\the\refno]\nref}
\def\nref#1{\xdef#1{[\the\refno]}\writedef{#1\leftbracket#1}%
\ifnum\refno=1\immediate\openout\rfile=refs.tmp\fi
\global\advance\refno by1\chardef\wfile=\rfile\immediate
\write\rfile{\noexpand\item{#1\ }\reflabeL{#1\hskip.31in}\pctsign}\findarg}
%	horrible hack to sidestep tex \write limitation
\def\findarg#1#{\begingroup\obeylines\newlinechar=`\^^M\pass@rg}
{\obeylines\gdef\pass@rg#1{\writ@line\relax #1^^M\hbox{}^^M}%
\gdef\writ@line#1^^M{\expandafter\toks0\expandafter{\striprel@x #1}%
\edef\next{\the\toks0}\ifx\next\em@rk\let\next=\endgroup\else\ifx\next\empty%
\else\immediate\write\wfile{\the\toks0}\fi\let\next=\writ@line\fi\next\relax}}
\def\striprel@x#1{} \def\em@rk{\hbox{}} 
\def\lref{\begingroup\obeylines\lr@f}
\def\lr@f#1#2{\gdef#1{\ref#1{#2}}\endgroup\unskip}

\def\addref#1{\immediate\write\rfile{\noexpand\item{}#1}} %now unnecessary
\def\startrefs#1{\immediate\openout\rfile=refs.tmp\refno=#1}
\def\xref{\expandafter\xr@f}\def\xr@f[#1]{#1}
\def\refs#1{\count255=1[\r@fs #1{\hbox{}}]}
\def\r@fs#1{\ifx\und@fined#1\message{reflabel \string#1 is undefined.}%
\nref#1{need to supply reference \string#1.}\fi%
\vphantom{\hphantom{#1}}\edef\next{#1}\ifx\next\em@rk\def\next{}%
\else\ifx\next#1\ifodd\count255\relax\xref#1\count255=0\fi%
\else#1\count255=1\fi\let\next=\r@fs\fi\next}
%

%
% this is ugly, but moore insists
\newwrite\ffile\global\newcount\figno \global\figno=1
\def\fig{fig.~\the\figno\nfig}
\def\nfig#1{\xdef#1{fig.~\the\figno}%
\writedef{#1\leftbracket fig.\noexpand~\the\figno}%
\ifnum\figno=1\immediate\openout\ffile=figs.tmp\fi\chardef\wfile=\ffile%
\immediate\write\ffile{\noexpand\medskip\noexpand\item{Fig.\ \the\figno. }
\reflabeL{#1\hskip.55in}\pctsign}\global\advance\figno by1\findarg}
\def\xfig{\expandafter\xf@g}\def\xf@g fig.\penalty\@M\ {}
\def\figs#1{figs.~\f@gs #1{\hbox{}}}
\def\f@gs#1{\edef\next{#1}\ifx\next\em@rk\def\next{}\else
\ifx\next#1\xfig #1\else#1\fi\let\next=\f@gs\fi\next}
\newwrite\lfile
{\escapechar-1\xdef\pctsign{\string\%}\xdef\leftbracket{\string\{}
\xdef\rightbracket{\string\}}\xdef\numbersign{\string\#}}

\def\writestop{\def\writestoppt{\immediate\write\lfile{\string\pageno%
\the\pageno\string\startrefs\leftbracket\the\refno\rightbracket%
\string\def\string\secsym\leftbracket\secsym\rightbracket%
\string\secno\the\secno\string\meqno\the\meqno}\immediate\closeout\lfile}}
\def\writestoppt{}\def\writedef#1{}
\def\seclab#1{\xdef #1{\the\secno}\writedef{#1\leftbracket#1}\wrlabeL{#1=#1}}
\def\subseclab#1{\xdef #1{\secsym\the\subsecno}%
\writedef{#1\leftbracket#1}\wrlabeL{#1=#1}}
\newwrite\tfile \def\writetoca#1{}
\def\leaderfill{\leaders\hbox to 1em{\hss.\hss}\hfill}
%	use this to write file with table of contents
\def\writetoc{\immediate\openout\tfile=toc.tmp 
   \def\writetoca##1{{\edef\next{\write\tfile{\noindent ##1 
   \string\leaderfill {\noexpand\number\pageno} \par}}\next}}}
%       and this lists table of contents on second pass

%
\catcode`\@=12 % at signs are no longer letters
%
%	Unpleasantness in calling in abstract and title fonts
\edef\tfontsize{\ifx\answ\bigans scaled\magstep3\else scaled\magstep4\fi}
\font\titlerm=cmr10 \tfontsize \font\titlerms=cmr7 \tfontsize
\font\titlermss=cmr5 \tfontsize \font\titlei=cmmi10 \tfontsize
\font\titleis=cmmi7 \tfontsize \font\titleiss=cmmi5 \tfontsize
\font\titlesy=cmsy10 \tfontsize \font\titlesys=cmsy7 \tfontsize
\font\titlesyss=cmsy5 \tfontsize \font\titleit=cmti10 \tfontsize
\skewchar\titlei='177 \skewchar\titleis='177 \skewchar\titleiss='177
\skewchar\titlesy='60 \skewchar\titlesys='60 \skewchar\titlesyss='60
\def\titlefont{\def\rm{\fam0\titlerm}% switch to title font
\textfont0=\titlerm \scriptfont0=\titlerms \scriptscriptfont0=\titlermss
\textfont1=\titlei \scriptfont1=\titleis \scriptscriptfont1=\titleiss
\textfont2=\titlesy \scriptfont2=\titlesys \scriptscriptfont2=\titlesyss
\textfont\itfam=\titleit \def\it{\fam\itfam\titleit}\rm}
 \ifx\answ\bigans\else scaled\magstep1\fi
\ifx\answ\bigans\def\abstractfont{\tenpoint}\else
\font\abssl=cmsl10 scaled \magstep1
\font\absrm=cmr10 scaled\magstep1 \font\absrms=cmr7 scaled\magstep1
\font\absrmss=cmr5 scaled\magstep1 \font\absi=cmmi10 scaled\magstep1
\font\absis=cmmi7 scaled\magstep1 \font\absiss=cmmi5 scaled\magstep1
\font\abssy=cmsy10 scaled\magstep1 \font\abssys=cmsy7 scaled\magstep1
\font\abssyss=cmsy5 scaled\magstep1 \font\absbf=cmbx10 scaled\magstep1
\skewchar\absi='177 \skewchar\absis='177 \skewchar\absiss='177
\skewchar\abssy='60 \skewchar\abssys='60 \skewchar\abssyss='60
\def\abstractfont{\def\rm{\fam0\absrm}% switch to abstract font
\textfont0=\absrm \scriptfont0=\absrms \scriptscriptfont0=\absrmss
\textfont1=\absi \scriptfont1=\absis \scriptscriptfont1=\absiss
\textfont2=\abssy \scriptfont2=\abssys \scriptscriptfont2=\abssyss
\textfont\itfam=\bigit \def\it{\fam\itfam\bigit}\def\footnotefont{\tenpoint}%
\textfont\slfam=\abssl \def\sl{\fam\slfam\abssl}%
\textfont\bffam=\absbf \def\bf{\fam\bffam\absbf}\rm}\fi
\def\tenpoint{\def\rm{\fam0\tenrm}% switch back to 10-point type
\textfont0=\tenrm \scriptfont0=\sevenrm \scriptscriptfont0=\fiverm
\textfont1=\teni  \scriptfont1=\seveni  \scriptscriptfont1=\fivei
\textfont2=\tensy \scriptfont2=\sevensy \scriptscriptfont2=\fivesy
\textfont\itfam=\tenit \def\it{\fam\itfam\tenit}\def\footnotefont{\ninepoint}%
\textfont\bffam=\tenbf \def\bf{\fam\bffam\tenbf}\def\sl{\fam\slfam\tensl}\rm}
\font\ninerm=cmr9 \font\sixrm=cmr6 \font\ninei=cmmi9 \font\sixi=cmmi6 
\font\ninesy=cmsy9 \font\sixsy=cmsy6 \font\ninebf=cmbx9 
\font\nineit=cmti9 \font\ninesl=cmsl9 \skewchar\ninei='177
\skewchar\sixi='177 \skewchar\ninesy='60 \skewchar\sixsy='60 
\def\ninepoint{\def\rm{\fam0\ninerm}% switch to footnote font
\textfont0=\ninerm \scriptfont0=\sixrm \scriptscriptfont0=\fiverm
\textfont1=\ninei \scriptfont1=\sixi \scriptscriptfont1=\fivei
\textfont2=\ninesy \scriptfont2=\sixsy \scriptscriptfont2=\fivesy
\textfont\itfam=\ninei \def\it{\fam\itfam\nineit}\def\sl{\fam\slfam\ninesl}%
\textfont\bffam=\ninebf \def\bf{\fam\bffam\ninebf}\rm} 
%
%---------------------------------------------------------------------
%

\hyphenation{anom-aly anom-alies coun-ter-term coun-ter-terms}
\def\inv{^{\raise.15ex\hbox{${\scriptscriptstyle -}$}\kern-.05em 1}}

\def\Dsl{\,\raise.15ex\hbox{/}\mkern-13.5mu D} %this one can be subscripted
\def\dsl{\raise.15ex\hbox{/}\kern-.57em\partial}

\font\bigit=cmti10 scaled \magstep1
 %pound sterling
\def\lspace{\ifx\answ\bigans{}\else\qquad\fi}
\def\lbspace{\ifx\answ\bigans{}\else\hskip-.2in\fi} % $$\lbspace...$$
\def\boxeqn#1{\vcenter{\vbox{\hrule\hbox{\vrule\kern3pt\vbox{\kern3pt
	\hbox{${\displaystyle #1}$}\kern3pt}\kern3pt\vrule}\hrule}}}
\def\mbox#1#2{\vcenter{\hrule \hbox{\vrule height#2in
		\kern#1in \vrule} \hrule}}  %e.g. \mbox{.1}{.1}
%	matters of taste
%\def\tilde{\widetilde} \def\bar{\overline} \def\hat{\widehat}
%
% some sample definitions
  %     curly letters
 \def\CC{{\cal C}}

\def\darr#1{\raise1.5ex\hbox{$\leftrightarrow$}\mkern-16.5mu #1}
 %pound sterling

 %puts a small half in a displayed eqn
\def\roughly#1{\raise.3ex\hbox{$#1$\kern-.75em\lower1ex\hbox{$\sim$}}}

\def\frac#1#2{{#1\over#2}}

\def\journal#1&#2(#3){\unskip, #1~\bf #2 \rm(19#3) }
\def\andjournal#1&#2(#3){\sl #1~\bf #2 \rm (19#3) }

\def\det{{\rm det}}

\catcode`\@=11\def\slash#1{\mathord{\mathpalette\c@ncel{#1}}}
\overfullrule=0pt
\def\steepslash{\c@ncel}
\def\frac#1#2{{#1\over #2}}

\def\:{\!:\!}
\def\inbar{\,\vrule height1.5ex width.4pt depth0pt}
\def\IQ{\relax\,\hbox{$\inbar\kern-.3em{\rm Q}$}}
\def\IB{\relax{\rm I\kern-.18em B}}
\def\IC{\relax\hbox{$\inbar\kern-.3em{\rm C}$}}
\def\IP{\relax{\rm I\kern-.18em P}}
\def\IR{\relax{\rm I\kern-.18em R}}
\def\ZZ{\relax\ifmmode\mathchoice
{\hbox{Z\kern-.4em Z}}{\hbox{Z\kern-.4em Z}}
{\lower.9pt\hbox{Z\kern-.4em Z}}
{\lower1.2pt\hbox{Z\kern-.4em Z}}\else{Z\kern-.4em Z}\fi}

\catcode`\@=12

%                      Zeitschriften:
\def\npb#1(#2)#3{{ Nucl. Phys. }{B#1} (#2) #3}
\def\plb#1(#2)#3{{ Phys. Lett. }{#1B} (#2) #3}
\def\pla#1(#2)#3{{ Phys. Lett. }{#1A} (#2) #3}
\def\prl#1(#2)#3{{ Phys. Rev. Lett. }{#1} (#2) #3}
\def\mpla#1(#2)#3{{ Mod. Phys. Lett. }{A#1} (#2) #3}
\def\ijmpa#1(#2)#3{{ Int. J. Mod. Phys. }{A#1} (#2) #3}
\def\cmp#1(#2)#3{{ Comm. Math. Phys. }{#1} (#2) #3}
\def\cqg#1(#2)#3{{ Class. Quantum Grav. }{#1} (#2) #3}
\def\jmp#1(#2)#3{{ J. Math. Phys. }{#1} (#2) #3}
\def\anp#1(#2)#3{{ Ann. Phys. }{#1} (#2) #3}
\def\prd#1(#2)#3{{ Phys. Rev. } {D{#1}} (#2) #3}
\def\ptp#1(#2)#3{{ Progr. Theor. Phys. }{#1} (#2) #3}
\def\aom#1(#2)#3{{ Ann. Math. }{#1} (#2) #3}

\def\br{\buildrel}

\def\C{{\bf C}}

\def\Z{{\bf Z}}
\def\P{{\bf P}}
\def\Q{{\bf Q}}
\def\Z{{\bf Z}}
\def\cA{{\cal A}}

\def\cF{{\cal F}}

\def\cI{{\cal I}}

\def\cL{{\cal L}}
\def\cM{{\cal M}}

\def\cO{{\cal O}}

\def\cR{{\cal R}}
\def\cS{{\cal S}}

\def\cX{{\cal X}}

\def\cicy#1(#2|#3)#4{\left(\matrix{#2}\right|\!\!
                     \left|\matrix{#3}\right)^{{#4}}_{#1}}

\def\ra{\rightarrow}

\def\bs{\bigskip}

\def\Box{{\,\lower0.9pt\vbox{\hrule 
\hbox{\vrule height 0.2 cm \hskip 0.2 cm  
\vrule height 0.2 cm}\hrule}\,}}

\global\newcount\thmno \global\thmno=0
\def\definition#1{\global\advance\thmno by1
\bigskip\noindent{\bf Definition \secsym\the\thmno. }{\it #1}
\par\nobreak\medskip\nobreak}
\def\question#1{\global\advance\thmno by1
\bigskip\noindent{\bf Question \secsym\the\thmno. }{\it #1}
\par\nobreak\medskip\nobreak}
\def\theorem#1{\global\advance\thmno by1
\bigskip\noindent{\bf Theorem \secsym\the\thmno. }{\it #1}
\par\nobreak\medskip\nobreak}
\def\proposition#1{\global\advance\thmno by1
\bigskip\noindent{\bf Proposition \secsym\the\thmno. }{\it #1}
\par\nobreak\medskip\nobreak}
\def\corollary#1{\global\advance\thmno by1
\bigskip\noindent{\bf Corollary \secsym\the\thmno. }{\it #1}
\par\nobreak\medskip\nobreak}
\def\lemma#1{\global\advance\thmno by1
\bigskip\noindent{\bf Lemma \secsym\the\thmno. }{\it #1}
\par\nobreak\medskip\nobreak}
\def\conjecture#1{\global\advance\thmno by1
\bigskip\noindent{\bf Conjecture \secsym\the\thmno. }{\it #1}
\par\nobreak\medskip\nobreak}
\def\exercise#1{\global\advance\thmno by1
\bigskip\noindent{\bf Exercise \secsym\the\thmno. }{\it #1}
\par\nobreak\medskip\nobreak}
\def\remark#1{\global\advance\thmno by1
\bigskip\noindent{\bf Remark \secsym\the\thmno. }{\it #1}
\par\nobreak\medskip\nobreak}
\def\problem#1{\global\advance\thmno by1
\bigskip\noindent{\bf Problem \secsym\the\thmno. }{\it #1}
\par\nobreak\medskip\nobreak}
\def\others#1#2{\global\advance\thmno by1
\bigskip\noindent{\bf #1 \secsym\the\thmno. }{\it #2}
\par\nobreak\medskip\nobreak}
\def\proof{\noindent Proof: }

\def\thmlab#1{\xdef #1{\secsym\the\thmno}\writedef{#1\leftbracket#1}\wrlabeL{#1=#1}}
%
% redefine \newsec so that all \thmno set to zero in a new section
%
\def\newsec#1{\global\advance\secno by1\message{(\the\secno. #1)}
\ifx\answ\bigans \vfill\eject \else \bigbreak\bigskip \fi  %if desired
\global\subsecno=0\thmno=0\eqnres@t\noindent{\bf\the\secno. #1}
\writetoca{{\secsym} {#1}}\par\nobreak\medskip\nobreak}
\def\eqnres@t{\xdef\secsym{\the\secno.}\global\meqno=1\bigbreak\bigskip}
\def\sequentialequations{\def\eqnres@t{\bigbreak}}\xdef\secsym{}
\nref\AM{P. Aspinwall and D. Morrison, {\it Topological
field theory and rational curves}, Commun. Math. Phys. 151
(1993), 245-262.}
\nref\AtiyahBott{M. Atiyah and R. Bott, {\it The moment map and
equivariant cohomology}, Topology 23 (1984) 1-28.} 
\nref\Batyrev{V. Batyrev, {\it Dual polyhedra and the mirror
symmetry for Calabi-Yau hypersurfaces in toric varieties},
Journ. Alg. Geom. 3 (1994) 495-535.}
\nref\BatyrevBorisov{V. Batyrev and L. Borisov, {\it 
On Calabi-Yau complete intersections in toric varieties},
alg-geom/9412017.}
\nref\BatyrevStraten{V. Batyrev and D. van Straten,
{\it Generalized hypergeometric functions and rational
curves on Calabi-Yau complete intersections in
toric varieties}, Comm. Math. Phys. 168 (1995) 495-533.}
\nref\BCKS{V. Batyrev, I. Ciocan-Fontanine, B. Kim
and D. van Straten, {\it Conifold transitions and mirror
symmetry for Calabi-Yau complete intersections in Grassmannian},
alg-geom/9710022.}
\nref\BCOV{M. Bershadsky, S. Cecotti, H. Ooguri, C. Vafa,
{\it Kodaira-Spencer theory of gravity and exact results
for quantum string amplitude},
Commun. Math. Phys. 165 (1994) 311-428.}
\nref\Bott{R. Bott, {\it A residue formula for holomorphic
vector fields}, J. Diff. Geom. 1 (1967) 311-330.}
\nref\CDGP{P. Candelas, X. de la Ossa, P. Green, and
L. Parkes, {\it A pair of Calabi-Yau manifolds as an
exactly soluble superconformal theory},
Nucl. Phys. B359 (1991) 21-74.}
\nref\CFKM{P. Candelas, X. de la Ossa, A. Font, S. Katz
and D. Morrison, {\it Mirror symmetry for two parameter 
models I}, hep-th/9308083.}
\nref\CaporasoHarris{L. Caporaso and J. Harris, {\it Parameter
spaces for curves on surfaces and enumeration of rational curves},
alg-geom/9608024.}
\nref\CrauderMiranda{B. Crauder and R. Miranda,
{\it Quantum cohomology of rational surfaces},
Progress in Mathematics 129, Birkhauser (1995) 33.}
\nref\DiFrancescoItzykson{P. DiFrancesco and C. Itzykson,
{\it Quantum intersection rings}, 
Progress in Mathematics 129, Birkhauser (1995) 81.}
\nref\Dijkgraaf{R. Dijkgraaf,
{\it Mirror symmetry and elliptic curves}, 
Progress in Mathematics 129, Birkhauser (1995) 149.}
\nref\Dubrovin{B. Dubrovin, {\it Geometry of 2D topological
field theory}, Springer LNM 1620 (1996) 120-348.}
\nref\Yau{{\it Essays on Mirror Manifolds},
ed. S.T. Yau, International Press.}
\nref\ES{G. Ellingsrud and S.A. Stromme, {\it The number
of twisted cubic curves on the general quintic threefolds},
in \Yau.}
\nref\FP{W. Fulton and R. Pandharipande, {\it Notes
on stable maps and quantum cohomology}, alg-geom/9608011.}
\nref\Getzler{E. Getzler, {\it The elliptic Gromov-Witten
invariants on $CP^3$}, alg-geom/9612009.}
\nref\GiventalII{A. Givental, {\it Equivariant Gromov-Witten
invariants}, alg-geom/9603021.}
\nref\GiventalIII{A. Givental, {\it A mirror theorem for
toric complete intersections}, alg-geom/9701016.}
\nref\GP{T. Graber and R. Pandharipande, {\it Localization of
Virtual classes}, alg-geom/9708001.} 
\nref\Ha{
R. Hartshorne, {\it Algebraic Geometry.}
 Graduate Texts in Mathematics
52.  Berlin, Heidelberg, New York
Springer  1977}
\nref\Hirzebruch{F. Hirzebruch, {\it Topological methods in
algebraic geometry}, Springer-Verlag, Berlin 1995, 3rd Ed.}
\nref\Hu{Y. Hu, {\it Moduli spaces of stable polygons and
symplectic structures on 
on $\bar\cM_{0,n}$}, Comp. Math. (to appear), alg-geom/9701011.}
\nref\Jinzenji{M. Jinzenji and M. Nagura, {\it Mirror symmetry
and exact calculation of $N-2$ point correlation function
on Calabi-Yau manifold embedded in $CP^{N-1}$}, hep-th/9409029.}
\nref\Kapranov{M. Kapranov, {\it Chow quotients of Grassmannian},
I.M. Gel'fand Seminar, AMS (1993).}
\nref\Kim{B. Kim, {\it On equivariant quantum cohomology},
IMRN 17 (1996) 841.}
\nref\KM{
F. Knudsen, and D. Mumford, {\it The projectivity of the moduli space
of stable curves  I},  preliminaries on `det' and `div'.   Math. Scand.
39 (1976) 19-55.}
\nref\MorrisonPlesser{D. Morrison and R. Plesser, {\it
Summing the instantons: quantum cohomology and mirror symmetry
in toric varieties}, alg-geom/9412236.}
\nref\HKTY{S. Hosono, A. Klemm,
           S. Theisen and S.T. Yau, {\it Mirror symmetry, mirror map
           and applications to complete intersection Calabi-Yau
           spaces}, hep-th/9406055.}
\nref\HLY{S. Hosono, B.H. Lian and S.T. Yau, {\it GKZ-generalized
hypergeometric systems and mirror symmetry of Calabi-Yau hypersurfaces},
alg-geom/9511001.}
\nref\Hsiang{W.Y. Hsiang, {\it On characteristic classes 
and the topological Schur lemma
from the topological transformation groups viewpoint},
Proc. Symp.
Pure Math. XXII. (1971) 105-112.}
\nref\SKatz{S. Katz, {\it On the finiteness of
rational curves on quintic threefolds}, Comp. Math.
60 (1986) 151-162.}
\nref\Kontsevich{M. Kontsevich, 
{\it Enumeration of rational curves via torus actions.}
In: The Moduli Space of Curves, ed. by
R\. Dijkgraaf, C\. Faber, G\. van der Geer, Progress in Math\.
vol\. 129, Birkh\"auser, 1995, 335--368.}
\nref\Li{J. Li, {\it Algebraic geometric interpretation
of Donaldson's polynomial invariants},
Journ. Diff. Geom. 37 (1993) 417-466.}
\nref\LibgoberTeitelboim
{A. Libgober and J. Teitelboim, Duke Math. Journ., Int.
Res.
Notices 1 (1993) 29.}
\nref\LLYII{B. Lian, K. Liu, and S.T. Yau, 
{\it Mirror Principle II}, in preparation.}
\nref\Manin{Yu.I. Manin, {\it Generating functions in algebraic geometry
and sums over trees.} In: The Moduli Space of Curves, ed. by
R\. Dijkgraaf, C\. Faber, G\. van der Geer, Progress in Math\.
vol\. 129, Birkh\"auser, 1995, 401--418.} 
\nref\ManinII{Yu.I. Manin, {\it Frobenius Manifolds, quantum cohomology,
and moduli spaces (Chapter I,II,III).}, Max-Planck Inst. preprint
MPI 96-113.}
\nref\Morrison{D. Morrison, {\it Picard-Fuchs Equations and
Mirror Maps for Hypersurfaces}, in \Yau.}
\nref\Okonek{ C. Okonek, M. Schneider and H. Spindler{\it Vector bundles on
complex projective spaces}, Progress in Math., Birkhauser.}
\nref\RuanTian{Y.B. Ruan and G. Tian, {\it A mathematical
theory of quantum cohomology}, Journ. Diff. Geom. Vol. 42,
No. 2 (1995) 259-367.}
\nref\Tian{G. Tian, {\it private communication}, November 1995.}
\nref\Voisin{C. Voisin, {\it A mathematical proof of
a formula of Aspinwall-Morrison},
Comp. Math. 104 (1996) 135-151.}
\nref\Witten{E. Witten, {\it Phases of N=2 theories in two
dimension}, hep-th/9301042.} 

\eject\vfill

\Title{}{Mirror Principle I}

\centerline{
Bong H. Lian,$^{1}$\footnote{}{$^1$~~Department of Mathematics,
Brandeis University, Waltham, MA 02154.}
Kefeng Liu,$^{2}$\footnote{}{$^2$~~Department of Mathematics,
Stanford University, Stanford, CA 94305.}
 and Shing-Tung Yau$^3$\footnote{}{$^3$~~Department of Mathematics,
Harvard University, Cambridge, MA 02138.} }

\vskip .2in

Abstract. 
We propose and study the following {\it Mirror Principle:}
certain sequences of multiplicative
equivariant characteristic classes
on Kontsevich's stable map moduli spaces can be computed
in terms of certain hypergeometric type classes.
As applications, we compute the equivariant Euler classes
of obstruction bundles 
induced by any concavex bundles -- including any
direct sum of line bundles -- on $\P^n$.
This includes proving the formula of
Candelas-de la Ossa-Green-Parkes hence
completing the program of Candelas et al, Kontesevich,
Manin, and Givental, to compute rigorously
the instanton prepotential function for the quintic in $\P^4$.
We derive, among many other examples,
the multiple cover formula for Gromov-Witten invariants
of $\P^1$, computed
earlier by Morrison-Aspinwall and by Manin in different approaches.
We also prove a formula for enumerating Euler classes which
arise in the so-called local mirror symmetry for some noncompact
Calabi-Yau manifolds.
At the end we  interprete an infinite dimensional transformation
group, called the mirror group, acting on Euler data,
as a certain {\it duality group} of the linear sigma model.

\Date{11/97} 
 
\baselineskip=20pt plus 2pt minus 2pt

\centerline{\bf Contents}\nobreak\medskip{\baselineskip=12pt
 \parskip=0pt\catcode`\@=11 \noindent {1.} {Introduction} \leaderfill{2} \par 
\noindent \quad{1.1.} {The Mirror Principle} \leaderfill{2} \par 
\noindent \quad{1.2.} {Enumerative problems and the Mirror Conjecture} \leaderfill{6} \par 
\noindent \quad{1.3.} {Acknowledgements} \leaderfill{8} \par 
\noindent {2.} {Euler Data} \leaderfill{9} \par 
\noindent \quad{2.1.} {Preliminaries and notations} \leaderfill{9} \par 
\noindent \quad{2.2.} {Eulerity} \leaderfill{11} \par 
\noindent \quad{2.3.} {Concavex bundles} \leaderfill{18} \par 
\noindent \quad{2.4.} {Linked Euler data} \leaderfill{25} \par 
\noindent \quad{2.5.} {The Lagrange map and mirror transformations} \leaderfill{28} \par 
\noindent {3.} {Applications} \leaderfill{35} \par 
\noindent \quad{3.1.} {The first convex example: The Mirror Conjecture} \leaderfill{38} \par 
\noindent \quad{3.2.} {First concave example: multiple-cover formula} \leaderfill{41} \par 
\noindent \quad{3.3.} {Second concave example: $K_{{\fam \bffam \tenbf P}^2}$} \leaderfill{42} \par 
\noindent \quad{3.4.} {A concavex bundle on ${\fam \bffam \tenbf P}^3$.} \leaderfill{43} \par 
\noindent \quad{3.5.} {A concavex bundle on ${\fam \bffam \tenbf P}^4$.} \leaderfill{44} \par 
\noindent \quad{3.6.} {General concavex bundles} \leaderfill{45} \par 
\noindent \quad{3.7.} {Equivariant total Chern class} \leaderfill{46} \par 
\noindent \quad{3.8.} {Concluding remarks} \leaderfill{47} \par 
 \catcode`\@=12 \bigbreak\bigskip} %input toc.tex
\writetoc %create toc.tmp

\newsec{Introduction}

\subsec{The Mirror Principle}

In section 2,  
we develop a general theory of Euler data,
and give many examples. In particular we introduce
the notions of a convex and a concave bundles on $\P^n$,
and show that they naturally give rise to Euler data.
In section 3 we apply our method to compute the equivariant
Euler classes, and their nonequivariant limits,
of obstruction bundles induced by a convex or a concave bundle.

We briefly outline our approach for computing multiplicative
equivariant characteristic classes
on stable map moduli. This outline also fixes some notations used later.
Our approach is partly inspired by Kontsevich's approach
using the torus action, Givental's idea of studying
equivariant Euler classes, Witten's idea of the linear sigma
model, and Candelas et al's idea of
mirror transformations. All four are syntheized with
what we call the {\it Mirror Principle}, which we now explain.

Let $M$ be a projective manifold and $\beta\in H^2(M,\Z)$.
Let $\bar\cM_{g,k}(\beta,M)$ be Kontsevich's stable map moduli space
of degree $\beta$, arithmetic genus $g$, with $k$ marked points
\Kontsevich.
For a good introduction to stable maps, see
the paper of Fulton-Pandharipande \FP.
Throughout this paper, we shall only deal with the case 
with $g=0$. 

We begin by analyzing two distinguished types of
fixed points under an induced torus $T$ action
on $\bar{\cal{M}}_{0,0}(d, \P^n)$. 
Both types of fixed points
reflect the structure of the stable map moduli
space. A smooth fixed point we consider
is a degree $d$ cover of a $T$-invariant $\P^1$ joining
two fixed points 
$p_i,p_j$ in $\P^n$. A singular fixed point
we consider is in the compactification divisor. It is given by gluing
together two 1-pointed maps
 $(f_1, C_1, x_1)\in \bar{\cal{M}}_{0,1}(r, \P^n)$ and 
$(f_1, C_2, x_2)\in \bar{\cal{M}}_{0,1}(d-r, \P^n)$ at the marked
points with $f_1(x_1)=f_2(x_2)=p_i\in \P^n$
($p_i$ being a $T$-fixed in $\P^n$), resulting in
a degree $d$ stable map $(f,C)$. We consider
two types of $T$-equivariant bundles $V$ on
$\P^n$, which we called convex and concave respectively
(Definition 2.7).
To be brief, we consider the convex case in this outline.
A convex bundle $V\ra\P^n$ induces on $\bar\cM_{0,0}(d,\P^n)$
an obstruction bundle $U_d$ whose fiber at $(f,C)$
is the section space $H^0(C,f^*V)$. First
we have the exact sequence over $C$
$$0\rightarrow f^*V\rightarrow f_1^*V\oplus f_2^*V\rightarrow
V|_{p_i}\rightarrow 0.$$
Passing to cohomology, we have
$$0\rightarrow H^0(C, f^*V)\rightarrow H^0(C_1, f_1^*V)\oplus  H^0(C_2, f_2^*V)
\rightarrow V|_{p_i}\rightarrow 0.$$
Hence we obtain a similar exact sequence
for the $U_d$ restricted to a suitable fixed point set.
Let $b$ be any
multiplicative equivariant characteristic class \Hirzebruch~
 for vector bundles.
The exact sequence on the fixed point set above gives
rise to the identity, which we call
the {\it gluing identity}:
$$b(V)\cdot b(U_d)= b(U_r)\cdot b(U_{d-r}).$$

Let $M_d:=\bar\cM_{0,0}((1,d),\P^1\times\P^n)$. This space
has a $G=S^1\times T$ action.
There is a natural equivariant contracting map $\pi:M_d\ra\bar{\cal
M}_{0,0}(d, \P^n)$ given by 
$$\pi:\ (f, C)\in M_d \mapsto (\pi_2\circ f, C')\in \bar{\cal
M}_{0,0}(d, \P^n)$$ where $\pi_2$ denotes the projection onto
the second factor of
$\P^1\times \P^n$. Here $C'=C$ if $(\pi_2\circ f, C)$ is still stable; and
if $(\pi_2\circ f, C)$ is unstable, this is the case when $C$ is of the form $C=C_1\cup
C_0\cup C_2$ with $\pi_2\circ f(C_0)$ a point in $\P^n$, then $C'$ is
obtained from $C$ by contracting the unstable component $C_0$.

Now via $\pi$, we pull back to $M_d$ all the information obtained above
on $\bar\cM_{0,0}(d,\P^n)$. The reason is
that there is an collapsing map 
$\varphi: M_d\rightarrow N_d=\P^{(n+1)d+n}$ 
which then allows us to perform computations on
the linear object $N_d$. We call {\it $N_d$ the linear sigma model
and $M_d$ the nonlinear sigma model}. There is a natural
$G$ action on $N_d$ such that $\varphi$ is $G$-equivariant. 
For example, to determine an equivariant cohomology class $\omega$ on
$N_d$, we only need to know its
restrictions $\iota^*_{p_{ir}}(\omega)$
to the $(n+1)(d+1)$ $G$-fixed points
$\{p_{ir}\}$ in $N_d$.
The corresponding weights of the $G$ action on $N_d$
are $\lambda_i+r\alpha$.
Let $Q_d$ be the push-forward of
$b(V_d)=\pi^* b(U_d)$ into $N_d$. Then
the classes $Q_d$ inherit the gluing identity
(Theorem 2.8)
$$\iota^*_{p_{i,0}}(Q_0)\cdot \iota^*_{p_{ir}}(Q_d)
=\iota^*_{p_{ir}}(Q_r)\cdot
\iota^*_{p_{ir}}(Q_{d-r})
$$
which is an identity in the ring $H^*_G(pt)=H^*(BG)$.
The sequence $Q_d$ is an example of what we call an {\it Euler
data} (Definition 2.3). We summarize the various ingredients 
used, now and later, in our constructions:
$$\matrix{  &  & V_d=\pi^*U_d & & U_d  & & \rho^*U_d  & & \cr
 &  & \downarrow          & & \downarrow  & & \downarrow & & \cr
         N_d&{\br\varphi\over\longleftarrow} & M_d & 
{\br\pi\over\longrightarrow} &
\bar\cM_{0,0}(d,\P^n) & {\br\rho\over\longleftarrow} & 
\bar\cM_{0,1}(d,\P^n) & {\br ev\over\longrightarrow} & \P^n}
$$
where $\rho$ forgets, and $ev$ evaluates at,
the marked point of a 1-pointed stable map. Also $U_d=\rho_!(ev^*V)$, the
push-forward of the pull-back of $V$.

The gluing identity is not enough to determine all $Q_d$.
In order to get further information
localize the $Q_d$ to a fixed point in $N_d$ whose inverse image in
$M_d$ is a {\it smooth} fixed point, 
and compute $Q_d$ at the special values of
$\alpha=(\lambda_i-\lambda_j)/d$. The property of $b(U_d)$ at a smooth
fixed point comes into play here. 
For simplicity in this outline,
let's consider for example the case when $b$ is the
equivariant Euler class $e_T$ and $V=\cO(l)$.
In this case
we find that, at $\alpha=(\lambda_i-\lambda_j)/d$, 
$$\iota^*_{p_{i,0}}(Q_d)
=\prod_{m=0}^{ld}
(l\lambda_i-m(\lambda_i-\lambda_j)/d),$$ 
which actually inherits
the same identity of $e_T(U_d)$ at a smooth fixed point in $\bar{\cal M}_{0, 0}(d, \P^n)$.  
This immediately tells us that the $Q_d$ should be compared with 
the sequence of classes
$$P:~~P_d=\prod_{m=0}^{ld}(l\kappa-m\alpha)$$ 
(Theorem 2.10), which has $\iota^*_{p_{i,0}}(P_d)=
\iota^*_{p_{i,0}}(Q_d)$ at
$\alpha=(\lambda_i-\lambda_j)/d$. 
$P$ is another example of an Euler data.
This Euler data  will then naturally
give rise to a generating function of hypergeometric type,
hence explains the very origin of these functions in enumerative
problems on stable moduli!
The same holds true for a large class of vector bundles $V$
on $\P^n$.

Finally, under a suitable bound
on $c_1(V)$, we can completely determine
the $Q_d$, 
in terms of $P_d$ by means of
a mirror transformation argument
(Theorem 3.9). We can in turn use $Q_d$ to compute
$e_T(U_d)$ and their nonequivariant limits (Theorem 3.2).
Our approach, thus, makes the roles of three objects and their
relationships quite transparent: certain fixed point sets,
equivariant multiplicative characteristic
classes, and series of hypergeometric type.
Our method works well for many characteristic classes
such as the Euler class and the total Chern class.
 
We summarize this {\it Mirror Principle}:
starting from an equivariant multiplicative characteristic class $b$, 
a vector bundle $V$ on $\P^n$ which induces a sequences of
bundles $U_d\ra\bar{\cal M}_{0,0}(d, \P^n)$,

\item{1.} (Euler Data) The behaviour of the $U_d$
 at a singular fixed point gives rise to the
gluing identity. In turn this defines an  Euler data $Q_d$
on the linear sigma model $N_d$.

\item{2.} (Linking) The behaviour of the $U_d$
at a smooth fixed point allows us to read off the restrictions
of $Q_d$ at certain fixed points of $N_d$ for special values of
weights. The restriction values determine a distinguished
Euler data $P_d$ to be compared with the $Q_d$.

\item{3.} (Mirror Transformations)
Compute the $Q_d$ and the $b(U_d)$ in terms of 
the $P_d$ explicitly by means of a mirror
transformation argument.

\bs

Our approach outlined here can
be applied to a rather broad range of cases by 
replacing $\P^n$ by other manifolds. They include
manifolds with torus action and their submanifolds. We will
study the cases of toric varieties, grassmannians, and
homogeneous manifolds in a forth-coming paper \LLYII.
On the other hand, we can even go beyond equivariant
multiplicative characteristic classes. In a future paper, 
we will study a sequence of equivariant classes $T_d$
of geometric origin
on stable moduli satisfying our gluing identity.

\subsec{Enumerative problems and the Mirror Conjecture}

For the remarkable history of the Mirror Conjecture,
see \Yau.
In 1990, Candelas, de la Ossa, Green, and Parkes conjectured a formula
for counting the number $n_d$ of rational curves 
in every  degree $d$ on a general
quintic in $\P^4$. Their computation is partly
inspired by an earlier construction of a mirror manifold
by Greene-Plesser. It has been conjectured earlier,
 by Clemens, that the
number of rational curves in every degree is finite.
The conjectured formula agrees with a classical result
in degree 1, an earlier computation
by S. Katz in degree 2, and has been verified in 
degree 3 by Ellingsrud-Stromme.  In 1994 following
some ideas of Gromov and Witten, 
Ruan-Tian introduced the notion of a symplectic
Gromov-Witten (GW) invariants.
Independently Kontsevich
proposed an algebraic geometric notion of 
GW invariants. Significant generalizations 
of his definition have been given by Kontsevich-Manin,
Li-Tian and Behrend-Fentachi.
A recent paper of
Li-Tian shows that the symplectic version
and the algebraic geometric version of the GW theory are
essentially the same in the projective catetory.
Beautiful applications of ideas from quantum cohomology 
have recently been done by Caporaso-Harris \CaporasoHarris,
Crauder-Miranda \CrauderMiranda, DiFrancesco-Itzykson 
\DiFrancescoItzykson~ and others,
solving many important enumerative  problems.
Significant connection between quantum cohomology
and the geometry of Frobenius manifolds appears in 
the work of Dubrovin \Dubrovin~ and Manin \ManinII.

Closer to mirror symmetry, a theorem of Manin says that 
the degree $k$ GW invariants for $\P^1$ (the so-called multi-cover
contribution) is given by $k^{-3}$. This was conjectured
by Candelas et al, and was  justified by
Aspinwall-Morrison using a different compactification.
See also  Voisin \Voisin.
According to Kontsevich, the number
$$K_d=\sum_{k|d} n_{d/k} k^{-3}$$
is the degree of the Euler class $c_{top}(U_d)$ for 
$U_d\ra \bar\cM_{0,0}(d,\P^4)$ induced by
$\cO(5)\ra\P^4$.
We shall call $K_d$ the degree $d$ Kontsevich-Manin number.

Using the torus action on $\P^4$ and 
the Atiyah-Bott localization formula, Kontsevich
computes the numbers
$K_d$ for $d=1,2,3,4$, and verifies that they
agree with the conjectured formula.
In some recent papers \GiventalII\GiventalIII,
 Givental has introduced
a number of beautiful ideas, among which is
an equivariant version of quantum cohomology theory
(see also Kim \Kim).
However, it is not clear how his approach gives a
a complete proof of the conjectured formula.   
More recently, Graber-Pandharipande
have also applied fixed point method to study Gromov-Witten
invariants of $\P^n$.

We believe that the machinery introduced in this paper
will be useful for many other enumerative problems,
aside from proving the formula of Candelas et al.
In fact we have  applied our machinery to problems
in local mirror symmetry proposed by C.Vafa, S. Katz, and others.
${}^4$\footnote{}
{S. Katz has informed us that A. Elezi has also
studied a similar problem.}

We now formulate one of our main theorems in this paper.
Let 
$$K_d=\int_{\bar\cM_{0,0}(d,\P^4)}c_{top}(U_d),~~~~
F(T)={5 T^3\over 6}+\sum_{d>0} K_d e^{dT}.$$
Consider the fourth order hypergeometric differential operator:
$$L:=({d\over dt})^4-5 e^t
(5{d\over dt}+1)\cdots
(5{d\over dt}+4).$$
By the Frobenius method, it is easy to show that
$$f_i:={1\over i!}({d\over dH})^i|_{H=0}
\sum_{d\geq0} e^{d(t+H)} 
{\prod_{m=1}^{5d}(5H+m)\over
\prod_{m=1}^d(H+m)^5},~~i=0,1,2,3,$$
form a basis of solutions to the differential equation $L\cdot f=0$.
Let
$$T={f_1\over f_0},~~~\cF(T)={5\over 2}({f_1\over f_0} {f_2\over f_0}
-{f_3\over f_0}).$$

\theorem{(The Mirror Conjecture) $F(T)=\cF(T)$.}
 
The transformation on the functions $f_i$ given by
the normalization
$$f_i\mapsto {f_i\over f_0}$$
and the change of variables
$$t\mapsto T(t)={f_1\over f_0}$$
are known as the {\it mirror transformation}. 
By the construction of Candelas et al, the functions
$f_0,..,f_3$ are periods of a family of Calabi-Yau threefolds.
By the theorem of Bogomolov-Tian-Todorov, these periods
in fact determine the complex structure of the threefold.
 
A similar Mirror Conjecture
formula holds true for a three dimensional
Calabi-Yau complete
intersection in a toric Fano manifold \LLYII. 
This will turn out to agree with the beautiful
construction of the mirror manifolds of Batyrev \Batyrev,
Batyrev-Borisov \BatyrevBorisov, as well as
the many mirror symmetry
computations of Morrison \Morrison, Libgober-Teitelboim
\LibgoberTeitelboim, Batyrev-van Straten
\BatyrevStraten, Candelas-Font-Katz-Morrison \CFKM,
Hosono-Klemm-Theisen-Yau \HKTY~ and
Hosono-Lian-Yau \HLY.

In this paper, to
make the ideas clear we restrict ourselves to the simplest case,
genus $0$ curves in some submanifolds of $\P^n$.
In a forth-coming paper \LLYII, we extend our discussions to
toric varieties, homogeneous manifolds,
their submanifolds, and for higher genus moduli spaces.
We hope to eventually understand from this point of view
the far reaching results for higher genus of 
Bershadsky-Cecotti-Ooguri-Vafa \BCOV, the beautiful
computations of Getzler \Getzler~ and Dijkgraaf \Dijkgraaf~
for elliptic GW invariants, and of
Batyrev-Ciocan Fontanine-Kim-van Straten \BCKS~
on Grassmannians.

\subsec{Acknowledgements}

We thank A. Todorov, A. Strominger, C. Vafa for helpful
discussions.
Our special thanks are due to J. Li who has been very
helpful throughout our project.

\newsec{Euler Data}

One of the key ingredients in our approach is
the linear sigma model, first introduced
by Witten \Witten, and later used to study mirror symmetry
by Morrison-Plesser \MorrisonPlesser, Jinzenji-Nagura \Jinzenji,
and others, resulting in
new insights into the origin of hypergeometric series.
In this paper, we consider the $S^1\times T$-equivariant
cohomology of the linear sigma model.

\subsec{Preliminaries and notations}

Let $T$ be an $r$-dimensional real torus with a
complex linear representation on $\C^{N+1}$. Let
$\beta_0,..,\beta_N$ be the weights of this action.
We consider the induced action of $T$ on $\P^N$,
and the $T$-equivariant cohomology with coefficients in $\Q$,
which we shall denote by $H^*_T(-)$.
Now $H^*_T(pt)$ is a polynomial algebra in $r$-variables,
and $\beta_i$ may be regarded as elements of
$H^2_T(pt)$.
{\it Throughout this paper, we shall follow
the convention that such generators have degree 1.} 
It is known that the equivariant cohomology of $\P^N$ is given by
\Hsiang
$$H_T(\P^N)=H_T(pt)[\zeta]/\left(\prod_{i=0}^N(\zeta-\beta_i)\right).$$
Here $\zeta$, which we shall call the equivariant hyperplane
class, is a fixed lifting of the hyperplane class of $\P^N$.
Each one-dimensional weight space in $\C^{N+1}$
becomes a fixed point $p_i$ in $\P^N$. {\it We shall identify
the rings $H^*_T(p_i)$ and $H_T^*(pt)=
H^*(BT)$.}  There are $N+1$ canonical
restriction maps $\iota^*_{p_i}:H_T(\P^N)\ra H_T(pt)$, given by
$\zeta\mapsto\beta_i$, $i=0,..,N$. There is also
a push-forward map $H_T(\P^N)\ra H_T(pt)$ given
by integration over $\P^N$. By the localization
formula, it is given by
$$\omega\mapsto
\int_{\P^N}\omega
= Res_\zeta{\omega\over\prod_{i=0}^N(\zeta-\beta_i)}.$$ 

Two situations
arise frequently in this paper.
First consider the standard action of $T=(S^1)^{n+1}$
on $\C^{n+1}$, and let $\lambda=(\lambda_0,..,\lambda_n)$ denote
the weights. 
On $\P^n$, there are $n+1$
isolated fixed points 
$p_0,..,p_n$. 
We shall denote the equivariant
hyperplane class by $p$, the canonical restrictions by 
$\iota^*_{p_i}:\omega
\mapsto\iota^*_{p_i}(\omega)=\omega(\lambda_i)$, 
and the push-forward  by
$$pf:H_T(\P^n)\ra H_T(pt)=\Q[\lambda].$$
We shall use the evaluation map $\lambda\mapsto 0$
on the ring $H^*_T(\P^n)$, and shall call this the
{\it nonequivariant limit}. In this limit, $p$ becomes
the ordinary hyperplane class $H\in H^*(\P^n)$.

We now consider the second situation. For each $d=0,1,2,..$,
consider the following
complex linear action of the group $G:=S^1\times T$
on $\C^{(n+1)(d+1)}$. We let the group act on
the $(ir)$-th coordinate line in $\C^{(n+1)(d+1)}$ by 
the weights $\lambda_i+r \alpha$,
$i=0,..,n$, $r=0,..,d$. Thus there are $(n+1)(d+1)$
isolated fixed points $p_{ir}$ on the projective space $\P^{(n+1)d+n}$, given
by those coordinate lines. In this case, we shall
denote the equivariant hyperplane class by $\kappa$,
the canonical restrictions by
$\iota^*_{p_{ir}}:\omega\mapsto\iota^*_{p_{ir}}(\omega)=
\omega(\lambda_i+r \alpha)$,
and the push-forward by 
$$pf_d:H_G(\P^{(n+1)d+n})\ra H_G(pt)=\Q[\alpha,\lambda].$$
Here we have abused the notation $\kappa$, using it to represent
a class in $H_G(\P^{(n+1)d+n})$ for every $d$. But it should present
no confusion in the context it arises.

Let $N_d$ be the space of nonzero $(n+1)$-tuple of degree $d$
homogeneous polynomials in two variables $w_0,w_1$, modulo scalar.
There is a canonical
way to identify $N_d$ with $\P^{(n+1)d+n}$. Namely,
a point $z\in\P^{(n+1)d+n}$ corresponds to the
polynomial tuple $[\sum_r z_{0r}w_0^r w_1^{d-r},..,
\sum_r z_{nr}w_0^r w_1^{d-r}]\in N_d$.
{\it This identification will be used throughout
this paper}. It is easy to
see that the natural $T$-action on $(n+1)$-tuples
together with the $S^1$-action on $[w_0,w_1]\in\P^1$ by weights
$(\alpha,0)$,
induces a $S^1\times T$-action on $N_d$ which
coincides with the $S^1\times T$-action on $\P^{(n+1)d+n}$
described earlier.

\definition{(Notations) 
We call the sequence of projective spaces $\{N_d\}$ 
the linear sigma model for $\P^n$. Here are
some frequently used notations: $G:=S^1\times T$,
$\cR:=\Q(\lambda)[\alpha]$, $\cR^{-1}:=\Q(\lambda,\alpha)$,
$\cR H^*_G(N_d):=H^*_G(N_d)\otimes_{\Q[\lambda,\alpha]}\cR$,
$\cR^{-1} H^*_G(N_d):=H^*_G(N_d)\otimes_{\Q[\lambda,\alpha]}\cR^{-1}$,
and $deg_\alpha~\omega$ means
the degree in $\alpha$ of $\omega\in\cR$.}

Obviously the maps $\iota^*_{p_{ir}}$, $pf_d$ defined linearly
over $\Q[\lambda,\alpha]$, 
can be extended $\cR$- or $\cR^{-1}$-linearly.
There are two natural equivariant maps between
the $N_d$, given by
$$\eqalign{
&I:N_{d-1}\ra N_d,~~~
[f_0,..,f_n]\mapsto[w_1f_0,..,w_1f_n]\cr 
&\bar{~~~}:N_d\ra N_d,~~~
[f_0(w_0,w_1),..,f_n(w_0,w_1)]\mapsto
[f_0(w_1,w_0),..,f_n(w_1,w_0)].}
$$
The second map induces 
on equivariant cohomology $\cR^{-1}H_G^*(N_d)$,
$\bar\kappa=\kappa-d\alpha$,
$\bar\alpha=-\alpha$, $\bar\lambda_i=\lambda_i$.
In particular
 any $x\in\cR$ has the form $x=x_- +x_+$ with $\bar x_\pm=\pm x_\pm$.
We also extend $\bar{~~~}$ to the power series ring $\cR[[e^t]]$
by leaving $t$ invariant.

Composing a chain $d$  $I$'s, we get a canonical
inclusion  $N_0=\P^n{\br I_d\over\longrightarrow} N_d$.
Note that the image of the fixed point $p_i$ is $p_{i,0}$.
For $\omega\in \cR^{-1} H_G^*(N_d)$, we shall denote
by $I^*_d(\omega)\in \cR^{-1} H_G^*(N_0)$ 
the restriction of $\omega$ to $N_0$.
Since $S^1$ acts trivially on $\P^n$, we can write
$$H^*_G(N_0)=H^*_T(\P^n)[\alpha].$$
In particular note that $H^*_T(\P^n)$ is invariant under $\bar{~~~}$,
and that $I_d^*(\kappa)=p$.

Obviously the set
of classes $\omega\in H^*_T(\P^n)$ with $\iota^*_{p_i}(\omega)\neq0$ for
all $i$, is closed under multiplication. We localize 
the ring $H^*_T(\P^n)$ by allowing to invert such elements $\omega$. 
We denote the resulting ring by
$H^*_T(\P^n)^{-1}$.

\definition{(Notations)
The degree in $\alpha$ of a class $\omega\in H^*_T(\P^n)^{-1}[\alpha]$
will be denoted by $deg_\alpha~\omega$.
A class $\Omega\in H^*_T(\P^n)^{-1}$ with
$\iota^*_{p_i}(\Omega)\neq0$ for all $i$ will be called invertible.
Throughout this paper, 
$\Omega$ will denote a fixed but arbitrary invertible class.
$\cS$ denotes the set of sequences of cohomology classes
$Q:~~Q_d 
\in\cR^{-1} H^*_G(N_d),
~~~d=1,2,...$
}

\subsec{Eulerity}

\definition{
We call a sequence  
$Q:~~Q_d\in \cR H^*_G(N_d),~~~d=1,2,...,$
an $\Omega$-Euler data if for all $d$, and $r=0,..,d$, $i=0,..,n$,
$$(*)~~~
\iota^*_{p_i}(\Omega)~\iota^*_{p_{i,r}}(Q_d)=
\overline{\iota^*_{p_{i,0}}(Q_r)}
~\iota^*_{p_{i,0}}(Q_{d-r}),$$
where $Q_0:=\Omega$.
We denote by $\cA^\Omega$ the set of $\Omega$-Euler data.
}
{\it When  dealing  only with
one fixed class $\Omega$ at a time,
we shall say Euler rather than $\Omega$-Euler,
and shall write $\cA$ for the set of Euler data.}

More explicitly, condition (*) can be written as
$\Omega(\lambda_i)Q_d(\lambda_i+r \alpha)=
\overline{Q_r(\lambda_i)}
Q_{d-r}(\lambda_i).$
Applying this at $r=d$, we find that
\eqn\dumb{ Q_d(\lambda_i+d \alpha)=
\overline{Q_d(\lambda_i)}.}
Putting $ \alpha=(\lambda_j-\lambda_i)/d$, we see that
$Q_d(\lambda_j)$~at~$\alpha=(\lambda_j-\lambda_i)/d$ coincides~with~
$Q_d(\lambda_i)$~at~$\alpha=(\lambda_i-\lambda_j)/d$.
Applying both (*) and \dumb~
at $\alpha=(\lambda_j-\lambda_i)/r$
(hence $\lambda_j=\lambda_i+r \alpha$), we get
$$\Omega(\lambda_i)
Q_d(\lambda_j)
=Q_r(\lambda_j)
Q_{d-r}(\lambda_i)~~at
~\alpha=(\lambda_j-\lambda_i)/r.
$$

\lemma{(Reciprocity Lemma)
If $Q$ is an Euler data, then
for $i,j=0,..,n$, $r=0,1,..,d$, $d=0,1,2,..$, we have
\item{(i)} $Q_d(\lambda_i+d \alpha)=
\overline{Q_d(\lambda_i)}$.
\item{(ii)} $Q_d(\lambda_j)$ at $\alpha=(\lambda_j-\lambda_i)/d$
coincides with
$Q_d(\lambda_i)$ at $\alpha=(\lambda_i-\lambda_j)/d$ for $d\neq0$.
\item{(iii)}
$\Omega(\lambda_i)
Q_d(\lambda_j)
=Q_r(\lambda_j)
Q_{d-r}(\lambda_i)$
at $\alpha=(\lambda_j-\lambda_i)/r$ for $r\neq0$.
}

{\it Example 1.} Let $l$ be a positive integer. Put 
$$
P:~~P_d=\prod_{m=0}^{ld}(l\kappa-m \alpha)\in H^*_G(N_d).
$$
It is straightforward to check that
$P$ is an $lp$-Euler data. 
We leave this as an exercise to the reader. This will arise naturally
in the problem of computing the equivariant Euler classes
of the obstruction bundles induced
by $\cO(l)\ra\P^n$. (See below.)

{\it Example 2.} Put $\Omega=p^{-2}$, and
$$
P:~~P_d=\prod_{m=1}^{d-1}(\kappa-m \alpha)^2\in H^*_G(N_d).
$$
This Euler data will arise in the problem of computing
the so-called multiple cover contributions, ie. the Gromov-Witten
invariants for $\P^1$.

{\it Example 3.} Put $\Omega=(-3p)^{-1}$, and
$$
P:~~P_d=\prod_{m=1}^{3d-1}(-3\kappa+m \alpha)\in H^*_G(N_d).
$$
This Euler data will arise in the problem of computing
the equivariant Euler classes
of the obstruction bundles induced
by the canonical bundle of $\P^2$.

{\it Example 4.} Put $\Omega=-1$, and
$$
P:~~P_d=\prod_{m=0}^{2d}(2\kappa-m \alpha)
\times\prod_{m=1}^{2d-1}(-2\kappa+m \alpha)\in H^*_G(N_d).
$$
This Euler data will arise in the problem of computing
the equivariant Euler classes
of the obstruction bundles induced
by $\cO(2)\oplus\cO(-2)$ on $\P^3$. 

{\it Example 5.} It is easy to show that if $Q$ is
$\Omega$-Euler and $a\in\Q(\lambda)$ is any nonzero element,
 then the data
 $Q':~~Q'_d=a\cdot Q_d$ is $a\Omega$-Euler.
Similarly, the data $Q':~~Q'_d=(-1)^{d+1}Q_d$ is
also $-\Omega$-Euler.

{\it Example 6.}
We observe that the set of Euler data
is a monoid, ie. it is closed under the product
$Q_dQ'_d$, and has the unit given by $Q_d=1$ for all $d$.
Hence the product of an $\Omega$-Euler with
an $\Omega'$-Euler data is an $\Omega\Omega'$-Euler
data.
In the geometrically setting,
this multiplicative property sometimes corresponds to taking
intersection of two suitable projective manifolds.
In this case, 
the class $\Omega\in H^*_T(\P^n)$ plays the role
of the equivariant Thom class of the normal bundle of
such a projective manifold.

{\it Example 7.} Let $Q_d=\kappa(\kappa-d\alpha)\in H_G^*(N_d)$.
Then it is again trivial to check that $Q$ is $\kappa^2$-Euler.

{\it Example 8.} Let $Q$ be an $\Omega$-Euler data,
and $Q'$ be  an $\Omega'$-Euler data. Suppose
$Q_d/Q_d'\in \cR H_G^*(N_d)$ for all $d>0$. Then
it is immediate that they form a sequence, denoted by $Q/Q'$,
which is $\Omega/\Omega'$-
Euler. As a special case,
let $Q_d=l^2\kappa(\kappa-d\alpha)$ as in Example 7,
and let $P$ as in Example 1. Then $P/Q$ is a 
$(lp)^{-1}$-Euler data. Example 2 is
is obtained by squaring this. 
Examples 3 and 4 can also be obtained
in a similar way.

{\it Example 9.} Introduce a formal variable $x$. We can extend
everything above by adjoining $x$, ie. by replacing the
ground field $\Q$ be the ring $\Q[x]$.
For example, it is easy to show that
$$
P:~~P_d=\prod_{m=0}^{ld}(x+l\kappa-m \alpha)\in H^*_G(N_d)[x]
$$
satisfies the gluing identity as in Example 1, thus is an Euler data in
a more general sense. Such Euler data will appear
in the computations of equivariant total Chern classes.

{\it Example 10.} 
Let $M_d^0:=\cM_{0,0}((1,d),\P^1\times\P^n)$ be the moduli
space of holomorphic maps $\P^1\ra\P^1\times\P^n$ of bidegree $(1,d)$. 
Recall that $M_d$ is the stable map compactification of $M_d^0$.
Each map $f\in M_d^0$ can be represented by 
$f:[w_0,w_1]\mapsto[w_1,w_0]\times[f_0(w_0,w_1),..,f_n(w_0,w_1)]$,
where $f_i$ are degree $d$ homogeneous polynomials. 
So there is an obvious map 
$$\varphi:M_d^0\ra N_d,~f\mapsto[f_0,..,f_n]$$
which is $G=S^1\times T$ equivariant.
For convenience
we {\it define $M_0:=N_0=\P^n$}.
With a bit of work (see below), it can be shown that the map $\varphi$
has an equivariant regular extension to $\varphi:M_d\ra N_d$.
Let $(f,C)\in M_d$. Then $C$ is an arithmetic genus 0 curve
of  the form
$C=C_0\cup C_1\cup\cdots \cup C_N$ such that $\pi_1\circ f:C_0{\br\sim
\over\ra}\P^1$, where
$\pi_1,\pi_2$ 
are projections from $\P^1\times \P^n$ to the first
and second factors. 
Each $C_j$, $j>0$, 
is glued to $C_0$ at some point $x_j\in C_0$.
The map $\pi_2\circ f: \ C_j\ra \P^n$ 
is of degree $d_j$ with $\sum_j d_j=d$, 
and $\pi_1\circ f: \ C_j\ra \P^1$ is constant map 
with image $\pi_1\circ f(x_j)\in \P^1$. 
If we denote by $[\sigma_0, \cdots, \sigma_n]$ the degree $d_0$
polynomials representing $\pi_2\circ f: \ C_0\ra \P^n$, then 
$\varphi: (f,C)
\mapsto [\sigma_0g, \cdots, \sigma_ng]$,
where $g=\prod_j(a_jw_0-b_jw_1)^{d_j}$ 
with $\pi_1\circ f(x_j)=[a_j, b_j]$.
Thus $\varphi$ collapses all but one component of $C$. 

The idea of using a collapsing map relating two moduli problems
is not new.
The map $\varphi$ was known to Tian \Tian, and a similar map
also appeared in \Li~
in which a collapsing map was used to relate 
two moduli spaces. 
The map $\varphi$ was also
used in \GiventalII. Similar maps have also
been studied in \Hu\Kapranov.

Let $V=\cO(l)$ for $l>0$, and consider the induced bundle
$U_d\ra\bar\cM_{0,0}(d,\P^n)$. Pulling
 this back via the projection $\pi$,
we get a bundle $V_d\ra M_d$ (see Introduction).
Let $\chi_d$ be the equivariant Euler class of $V_d$.
We can now push-forward these classes for $d=1,2,..$
via the equivariant map $\varphi$ and obtain
a sequence
$\varphi_!(\chi_d)$. The following theorem will be
a special case of a general theorem proved in the next subsection.

\theorem{The sequence $\varphi_!(\chi_d)\in H^*_G(N_d)$ above
is an $lp$-Euler data.}
\thmlab\FirstExampleAdm

We now return to the map $\varphi$. The reader who wishes to
skip technical details can safely
omit the proof. 

\lemma{
There exists a morphism $\varphi:\ M_d\rightarrow N_d$.
Moreover $\varphi$ is equivariant with respect to the
induced action of $S^1\times T$.}
\thmlab\Regularity
\def\Pf{{\bf P}^n}
\def\Po{{\bf P}^1}
\let\po=\Po
\let\pf=\Pf
\def\lra{\longrightarrow}
\def\cL{{\cal L}}
\def\cO{{\cal O}}
\def\cF{{\cal F}}
\def\cS{{\cal S}}

\def\cX{{\cal X}}
\def\sta{^{\ast}}
\def\mor{{\rm Mor}\,}
\def\mh{\!:\!}

\def\cM{{\cal M}}
\def\upmo{^{-1}}
\def\tor{{\rm Tor}\,}
\def\CC{{\bf C}}

\proof The following proof is given by J. Li.
Let $M_d $ be the moduli space of stable morphisms $f: C\to \Po\times\Pf$ from
arithmetic genus 0 curves to $\Po\times \Pf$ of bi-degree $(1,d)$, and
let $N_d$ be the space of equivalence classes of
$(n+1)$-tuples $(f_0,\ldots,f_n)$, where
$f_i$ are degree $d$ homogeneous polynomials in two variables,
and $(f_0,\ldots,f_n)\sim(f_0',\ldots,f_n')$ if there is a
constant $c\ne0$ such that $f_i=c\cdot f_i'$ for all $i$.
We first define the morphism $\varphi: M_d \to N_d$. For convenience, we let
$\cS$ be the category of all schemes of finite type (over $\bf C$)
and let
$$\cF: \cS\lra ({\rm Set})
$$
be the the contra-variant functor that send any $S\in \cS$ to the
set of families of stable morphisms
$$F: \cX\lra \po\times\pf\times S
$$
over $S$, where $\cX$ are families of connected arithmetic
genus 0 curves, modulo the obvious equivalence relation. Note that
$\cF$ is represented by the moduli stack $M_d $. Hence to define $\varphi$
it suffices to define a transformation
$$\Psi: \cF\lra \mor(-,N_d).
$$
We now define such a transformation. Let $S\in\cS$ and let
$\xi\in \cF(S)$ be represented by
$(\cX,F)$. We let $p_i$ be the
composite of $F$ with the $i$-th projection of
$\po\times\pf\times S$ and let $p_{ij}$ be the composite of $F$ with
the projection from $\po\times\pf\times S$ to the product of its
$i$-th and $j$-th components.
We consider the sheaf $p_2\sta\cO_{\pf}(1)$ on $\cX$ and
its direct image sheaf
$$\cL_{\xi}=p_{13\ast} p_2\sta\cO_{\pf}(1).
$$
We claim that $\cL_{\xi}$ is flat over $S$. Indeed, by argument in
the proof of Theorem 9.9 in \Ha, it suffices to show that
$\pi_{S\ast}(\cL_{\xi}\otimes \pi_{\po}\sta\cO_{\po}(m))$ are locally free
sheaves of $\cO_S$-modules for $m\gg 0$,
where $\pi_{\po}$ and $\pi_S$ are the first and the
second projections of $\po\times S$. Clearly,
this sheaf is isomorphic to
$p_{3\ast}(p_2\sta\cO_{\pf}(1)\otimes p_1\sta\cO_{\po}(m))$,
which is locally free because
$$R^ip_{3\ast} (p_2\sta\cO_{\pf}(1)\otimes p_1\sta\cO_{\po}(m))=0
$$
for $i>0$ and $m\gg0$. For the same reasoning, the sheaves $\cL_{\xi}$
satisfy the following
base change property: let $\rho: T\to S$ be any base change and
let $\rho\sta(\xi)\in\cF(T)$ be the pull back of $\xi$. Then there is a
canonical isomorphism of sheaves of $\cO_T$-modules
\eqn\eqTwo{
 \cL_{\rho\sta(\xi)}\cong ({\bf 1}_{\po}\times \rho)\sta
\cL_{\xi}.
}

Since $\cL_{\xi}$ is flat over $S$, we can
define the determinant line bundle of $\cL_{\xi}$, denoted by $\det(\cL_{\xi})$
\KM.
The sheaf $\det(\cL_{\xi})$ is an invertible sheaf over $\po\times S$.
Using the Riemann-Roch theorem, one computes that its
degree along fibers over $S$ are
$d$. Further, because $\cL_{\xi}$ has rank one, there is a canonical
homomorphism
\eqn\eqThree{
\cL_{\xi}\lra\det(\cL_{\xi}),
}
so that its kernel is the torsion subsheaf of $\cL_{\xi}$,
denoted by $\tor(\cL_{\xi})$. Let $w_0,\ldots,w_n$ be the homogeneous
coordinate of $\pf$ chosen before. $w_0,\ldots,w_n$ form a
basis of $H^0(\pf,\cO_{\pf}(1))$. Then
their pull backs provide a collection of canonical sections of
$\cL_{\xi}$, and hence a collection of canonical sections
$$\sigma_{\xi,0},\ldots,\sigma_{\xi,n} \in H^0(S,\pi_{S\ast}\det(\cL_{\xi})).
$$
based on \eqThree.
Then after fixing an isomorphism
\eqn\eqOne{
\det(\cL_{\xi})\cong \pi_S\sta \cM\otimes\pi_{\po}\sta\cO_{\po}(d)
}
for some invertible sheaf $\cM$ of $\cO_S$-modules,
we obtain a section of
$$ \pi_{S\ast}(\pi_{\po}\sta\cO_{\po}(d))\otimes_{\cO_S}\cM\equiv
H^0_{\po}(\cO_{\po}(d))\otimes_{{\bf C}}\cM.
$$
Finally, we let $w_0, w_1$ be the homogeneous coordinate of $\po$.
Then the space $H^0_{\po}(\cO_{\po}(d))$ is the space of degree $d$
homogeneous polynomials in variables $w_0$ and $w_1$. This way, we obtain
a morphism
$$\Psi(S): S\lra N_d
$$
that is independent of the isomorphisms \eqOne.
It follows from the base change property \eqTwo~ that the collection
$\Psi(S)$ defines a transformation
$${\bf \Psi}: \cF\lra \mor(-,N_d),
$$
thus defines a morphism $\varphi$ as desired.

It remains to check that for any $w\in S$, the sections
$$\sigma_{\xi,0}(w),\ldots, \sigma_{\xi,n}(w)\in H^0(\po, \det(\cL_{\xi})
\otimes_{\cO_S} k_w)
$$
has the described vanishing property. Because of the base change property
of $\cL_{\xi}$, it suffices to check the statement when $S$ is
a point and $\xi\in\cF_d(S)$ is the stable map
$f\mh C\to \po\times\pf$. Let $x_1,\ldots,x_N$ be
the set of points in $\po$ so that
$p_1\mh C\to\po$, where $p_1=\pi_{\po}\circ f$,
is not flat over these points. Let $C_i$ be
$p_1\upmo(x_i)$ and
let $m_i$ be the degree of $f([C_i])\in H_2(\pf)$. Then
$\cL_{\xi}=p_{1\ast}p_2\sta\cO_{\pf}(1)$ is locally free away from
$x_1,\ldots,x_N$ and has torsion of length $m_i$ at $x_i$.
Then $\cL_{\xi}/\tor(\cL_{\xi})$ is locally free of degree
$k-\sum m_i$. It is direct to check that
the canonical inclusion
$$ \cL_{\xi}/\tor(\cL_{\xi})\lra \det(\cL_{\xi})\cong\cO_{\po}(d)
$$
has cokernel supported on the union of $x_1,\ldots,x_N$ whose
length at $x_i$ is exactly $m_i$. The statement about the vanishing of
$\sigma_{\xi,0}(w),\ldots, \sigma_{\xi,n}(w)$ follows immediately.

The fact that $\varphi: M_d \to N_d$ is $(\CC\sta)^{n+1}\times(\CC\sta)^2$-equivariant
as stated follows immediately from the fact that $\varphi$ is induced by
the transformation $\Psi$ of functors.
This completes the proof. $\Box$

\subsec{Concavex bundles}

\definition{We call a $T$-equivariant vector bundle $V\ra\P^n$ convex (resp.
concave) if the $T$-equivariant Euler class $e_T(V)$ is invertible
and if
$H^1(C,f^*V)=0$ (resp. $H^0(C,f^*V)=0$)
for every 0-pointed genus 0 stable map $f:C\ra\P^n$. We call $V$
a concavex bundle if it is a direct sum of a convex and a concave bundles.
We denote by $V^\pm$ the convex and concave summands of $V$.
By convention, we consider the zero bundle to be both convex and concave
so that concavexity includes both convexity and concavity.}

The convexity of a bundle is analogous to the
notion of convexity of a projective manifold introduced by
Behrend-Manin.

For example $\cO(l)\ra\P^n$ is convex if $l>0$, and concave if $l<0$.
Given any concavex
bundle $V\ra\P^n$, we have a sequence of induced
bundles 
$$U_d\ra\bar\cM_{0,0}(d,\P^n)$$
whose fiber at $(f,C)$
is the  space $H^0(C,f^*V^+)\oplus H^1(C,f^*V^-)$.
Pulling back $U_d$ via the contracting map $\pi:M_d\ra
\bar\cM_{0,0}(d,\P^n)$, we get a sequence
of bundles
$$V_d:=\pi^*U_d\ra M_d.$$
We denote by $\chi^V_d$ the equivariant Euler class of $V_d$,
and by $e_T(V)$ the equivariant Euler class of $V$.
We also introduce the notations:
$$\eqalign{
\Omega^V&:=
{e_T(V^+)\over e_T(V^-)}\cr
Q_d&:=\varphi_!(\chi_d^V), ~~Q_0:=\Omega^V.}
$$
By {\it convention},
if $V$ is the zero bundle, 
we set $\chi^V_d=1$, $e_T(V)=1$, $\Omega^V=1$.

\theorem{The sequence $\varphi_!(\chi^V_d)\in H_G^*(N_d)$ is an
$\Omega^V$-Euler data.}
\thmlab\GeneralEulerTheorem
\proof
We first discuss some preliminaries.
Let $M$ and $N$ be two compact smooth manifolds with the action of a torus
$T$, and $\varphi:\ M\rightarrow N$ be an equivariant map. Let $F$ be one
component of the fixed submanifold in $N$ and $i_F$ be the inclusion map
$F$ in $N$. Let $\phi_F={i_F}_!(1)\in H^*_T(N)$ denote the equivariant Thom class of
the normal bundle of $F$ in $N$. We then have, for any $\omega\in  H^*_T(M)$ 

$$\int_M\omega \varphi^*(\phi_F)=\int_N\varphi_!(\omega)\phi_F=\int_{F}i_F^*(\varphi_!(\omega)).$$

On the other hand, let $\{ P\}$ be the components of the fixed submanifold
contained in $\varphi^{-1}(F)$.
By  the Atiyah-Bott \Bott\AtiyahBott~ localization formula on $M$, we get 

$$\sum_P\int_{P}i_{P}^*(\omega \varphi^*(\phi_F))/e_T(P/M)=\int_{F}i_F^*(\varphi_!(\omega)).$$Here $e_T(P/M)$ denotes the
equivariant Euler class of the normal bundle of $P$ in $M$.
The reason is that the contribution of the fixed point sets not contained
in $\varphi^{-1}(F)$ is clearly zero. Actually assume $Q$ is a component
not contained inside $\varphi^{-1}(F)$, its contribution to the
localization is given by 
$$\int_{Q}i_{Q}^*(\omega\varphi^*(\phi_F))/e_T(Q/M).$$But by the naturality
of the pull-backs, we have 
 $$i^*_Q\varphi^*(\phi_F)=\varphi_0^*i_E^*(\phi_F)=0$$ where $E=\varphi(Q)$
is a fixed submanifold in $N$ and $\varphi_0$ denotes the restriction of
$\varphi$ to $Q$. 
Note that if $F$ is an isolated point, then $i^*_P\varphi^*(\phi_F)$ can be
pulled out of the integral.

The above formula will be applied to the collapsing map $\varphi:\
M_d\rightarrow N_d$.
All manifolds involved here are at worst orbifolds with finite quotient singularities, so the localization formula remains
valid without any change as long as we consider the corresponding integrals in the orbifold sense.

We consider the $S^1$-action on the 
$\P^1$ factor in $\P^1\times \P^n$ with weights $\alpha, 0$.
Combining with the natural $T$-action on $\P^n$, we get the naturally 
induced $G=S^1\times T$-actions on $M_d$ and $N_d$, with respect to which
the collapsing map $\varphi$ is equivariant. As described in section 2.1, the $G$-fixed points in $N_d$ are all 
of the form 
$$p_{ir}=[0, \cdots, 0, w_0^{r}w_1^{d-r}, 0,\cdots, 0]$$ 
in which the only nonzero term is in the $i$-th position. 

For each $r> 0$, let $\{ F_r\}\subset \bar\cM_{0,1}(r,\P^n)$ denote the
$T$-fixed point components in $\bar\cM_{0,1}(r,\P^n)$ with the marked point mapped
onto the fixed point $p_i$ in $\P^n$. Let $N(F_r)=N_{F_r/\bar\cM_{0,1}(r,\P^n)}$ denote the normal bundle of $F_r$ in $\bar\cM_{0,1}(r,\P^n)$.

Let $\pi_1,\pi_2$ be the projections 
from $\P^1\times\P^n$
onto the first and
second factors. 
From the construction of 
$\varphi$, we see that the $G$-invariant submanifold that is mapped to $p_{ir}$ consists 
of the following degree $(1, d)$ stable maps $f:\ C\rightarrow \P^1\times \P^n$ with 
$C=C_1\cup C_0\cup C_2$. Here $C_0\simeq \P^1$ and 
 $$\pi_2\circ f (C_0)=[0, \cdots, 0,1,0\cdots, 0]=p_i\in \P^n$$
where $1$ is at the $i$-th position. The map $\pi_1\circ f:
\ C_0\rightarrow \P^1$ is an isomorphism and maps $x_1=C_0\cap C_1$ 
and $x_2=C_0\cap C_2$ to $0$ and $\infty$ respectively. Actually 
$$\pi_1\circ f(C_1)=0, ~~~ \pi_1\circ f(C_2)=\infty ~~in~\P^1,$$ 
ie. the curves $C_1$ and $C_2$ are respectively
 mapped to the points $0$ and $\infty$ of $\P^1$.

The maps $\pi_2\circ f$ restricted 
to $C_j$ for $j=1, 2$ are stable maps in $\bar{{\cal M}}_{0,1}(r, \P^n)$
and $\bar{{\cal M}}_{0,1}(d-r, \P^n)$ respectively. We consider $F_r\times
F_{d-r}$ as a $G$-fixed submanifold of $M_d$ by gluing each pair 
to $C_0$ at $x_1$ and $x_2$ respectively as above. It is easy to see that
$\{F_r\times F_{d-r}\}$ are the $G$-fixed point sets in 
$M_d$ whose image under $\varphi$ is the fixed point $p_{ir}$. 

We first consider a convex bundle $V$ on $\P^n$,
and the case $r\neq 0,
d$. Then we have  
\eqn\Qir{
Q_d(\lambda_i+r\alpha)=\int_{N_d}\phi_{p_{ir}}Q_d=    \int_{M_d}\varphi^*(\phi_{p_{ir}})\chi^V_d.}
 Here $\phi_{p_{ir}}$ denote the equivariant Thom class of the $G$-fixed
 point $p_{ir}$ in $N_d$. We will apply the
 localization formula to compute the 
 right hand side of \Qir. First we need to know the normal
 bundle of the fixed points, which is, in the equivariant $K$-group of 
$F_r\times F_{d-r}$ \Kontsevich, \GP,  
$$\eqalign{
N_{F_r\times F_{d-r}/M_d}
&=N(F_r)+N(F_{d-r})+[H^0(C_0, (\pi_1\circ f)^*T\P^1)]\cr
&+ [L_r\otimes T_{x_1}C_0]+[L_{d-r}\otimes T_{x_2}C_0]-[T_{p_i}P^n]-[A_{C_0}]
}$$
where $L_r, L_{d-r}$ are the line bundles on 
$\bar\cM_{0,1}(r,\P^n), ~ \bar\cM_{0,1}(d-r,\P^n)$
 respectively, whose fiber at a stable map is the 
tangent line at the corresponding marked points $x_1$ and $x_2$. They
correspond to the deformation of the nodal points $x_1$ and
$x_2$. The term $H^0(C_0, (\pi_1\circ f)^*T\P^1)$ corresponds to the
deformation of $f$ restricted to $C_0$. The term $[A_{C_0}]$ is the bundle representing 
the infinitesimal automorphism of $C_0$ fixing the two 
points $x_1, x_2$. The term $-[T_{p_i}P^n]$ comes from gluing $F_r$ and
$F_{d-r}$ onto $C_0$ and the property that $\pi_2\circ f(C_0)=p_i$. 

This gives the following formula for the corresponding equivariant Euler classes:
$$\eqalign{
 e_T(F_r\times F_{d-r}/M_d)&=e_T(N(F_r))e_T(N(F_{d-r}))e_T(L_r\otimes T_{x_1}C_0)e_T(L_{d-r}\otimes T_{x_2}C_0)\cr
&\times e_T(T_{p_i}\P^n)^{-1}e_T(H^0(C_0, (\pi_1\circ f)^*T\P^1))e_T^{-1}(A_{C_0}).}
$$
Each term in this formula can be explicitly calculated. We clearly have $ e_T(T_{p_i}\P^n)= \prod_{j\neq i}(\lambda_i-\lambda_j)$; the weights of $T_{x_1}C_0$ and $ T_{x_2}C_0$ are $\alpha$ and $-\alpha$ respectively, therefore

$$e_T(L_r\otimes T_{x_1}C_0)=\alpha+c_1(L_r),~~~ e_T(L_{d-r}\otimes T_{x_2}C_0)=-\alpha+c_1(L_{d-r})$$
where $c_1(L_r), ~ c_1(L_{d-r})$ are the restriction to $F_r$ and $F_{d-r}$ of 
the equivariant Chern classes of the line bundles $L_r$ and $L_{d-r}$ with respect to the induced $T$ 
actions on 
$\bar\cM_{0,1}(r,\P^n)$ and $ \bar\cM_{0,1}(d-r,\P^n)$. To compute $e_T(H^0(C_0, (\pi_1\circ f)^*T\P^1))$ and $e_T(A_{C_0})$, first note that we have the standard exact sequence

$$0\rightarrow O\rightarrow \cO(1)\otimes \C^2\rightarrow T\P^1\rightarrow
0,$$ 
with $O$ being the trivial bundle. From this we get 

$$0 \rightarrow O\rightarrow H^0(C_0, \cO(1))\otimes \C^2\rightarrow H^0(C_0,(\pi_1\circ f)^*(T\P^1)) \rightarrow 0.$$

The weights of $ H^0(C_0, \cO(1))$ with basis $\{ w_0,\, w_1\}$ are $ \alpha
, \ 0$,
the weights of $\C^2$ with basis $\{ \frac{\partial}{\partial w_0}, \,
\frac{\partial}{\partial w_1}\}$ are $- \alpha, 0$ and the weight of $O$ is $0$. Therefore one finds that the weights of $H^0(C_0,(\pi_1\circ f)^*(T\P^1))$ are $ \alpha,\ - \alpha,\ 0$.

For $[A_{C_0}]$, we have the exact sequence 

$$0 \rightarrow A_{C_0}\rightarrow H^0(C_0,(\pi_1
\circ f)^*(T\P^1))\rightarrow T_{x_1}C_0\oplus T_{x_2}C_0\rightarrow 0.$$ The weights of $T_{x_1}C_0$ and $T_{x_2}C_0$ are $ \alpha$ and $- \alpha$ respectively. So $[A_{C_0}]$ contributes a $0$ weight space which cancels with the $0$ weight space of $[H^0(C_0,(\pi_1\circ f)^*(T\P^1))]$. Therefore we will ignore the zero weights in our formulas and write as 
$$e_T(H^0(C_0, (\pi_1\circ f)^*(T\P^1))e_T^{-1}(A_{C_0})= - \alpha^2.$$ 

When $V_d$ is restricted to $F_r\times F_{d-r}$ considered as
a fixed point set of $M_d$ as before, we have the exact sequence:

$$0\rightarrow V_d\rightarrow V_{r}|_{F_r}\oplus V_{d-r}|_{F_{d-r}}\rightarrow
V|_{p_i}\rightarrow 0.$$ Note that $V_{r}|_{F_r}$ and $V_{d-r}|_{F_{d-r}}$
is the same as $\rho^*U_{r}|_{F_r}$ and $\rho^*U_{d-r}|_{F_{d-r}}$ which respectively are the restrictions
to $F_r$ and $F_{d-r}$ of the pull-backs to $\bar{\cal M}_{0,
1}(r, \P^n)$ and $\bar{\cal M}_{0, 1}(d-r , \P^n)$ of the corresponding
bundles on $\bar\cM_{0,0}(r,\P^n)$ and $\bar\cM_{0,0}(d-r,\P^n)$. Here
$V|_{p_i}$ denotes the fiber of $V$ at $p_i\in \P^n$.

Here comes the important point. The multiplicativity of equivariant Euler classes gives us 
$$\Omega^V(\lambda_i)\cdot \chi_d^V=\chi^V_{r}\cdot \chi^V_{d-r}=
\rho^* e_T(U_r)\cdot \rho^* e_T(U_{d-r}),$$
when restricted to 
$F_r\times F_{d-r}$. Here $\rho:\bar\cM_{0,1}(d,\P^n)\ra
\bar\cM_{0,0}(d,\P^n)$ is the forgetting map
(same notation for all $d$).
Note that the above equality is just the pull-back via $\pi$ from $\bar{\cal
M}_{0,0}(d, \P^n)$ of the gluing identity discussed in the Introduction. 

For the case of $r=0$ or $d$, there is only one of the 
curves $C_1$ or $C_2$, that is $C$ is of the 
form $C_0\cup C_2$ or $C_1\cup C_0$. In this case we identify $F_d$ as the fixed point set in $M_d$ by gluing its marked point to $C_0$ at $x_1$ or $x_2$. The normal bundle in these two cases are respectively given by 

$$N_{F_d/M_d}=N(F_d)+[H^0(C_0, (\pi_1\circ f)^*T\P^1)]+[L_d\otimes T_{x_j}
C_0]-[A_{C_0}]$$ in the $K$-group of $F_d$. Here $L_d$ is the restriction to
$F_d$ of the line bundle on $\bar\cM_{0,1}(d,\P^n)$ whose fiber is the
tangent line at the marked point. For simplicity we write $L$ as $L_d$ in
the following. For $j=1, \, 2$, $T_{x_j}C_0$ is the
tangent line of $C_0$ at the corresponding marked point $x_j$. In these two
cases, one easily shows in the same way as above that the 
term $e_T(A_{C_0})$, except the $0$ weight, contributes one nonzero weight $- \alpha$ or $ \alpha$ respectively. Its $0$ weight space still cancels with the $0$ weight space of $[H^0(C_0, (\pi_1\circ f)^*T\P^1)]$.

By putting all of the above computations together and combining
with \Qir, we get, for $r\neq 0,\, d$, 
\eqn\dumbI{\eqalign{
\Omega^V(\lambda_i)Q(\lambda_i+r\alpha)
&=\Omega^V(\lambda_i) \int_{M_d}\varphi^*(\phi_{p_{ir}})\chi_d\cr
&=- \alpha^{-2}\prod_{j\neq i}
(\lambda_i-\lambda_j)e_T(p_{ir}/N_d)~~
\sum_{F_r}\int_{{F_r}}\frac{\rho^* e_T(U_r)}{e_T(N(F_r))( \alpha+c_1(L_r))}\cr
&\times\sum_{F_{d-r}}\int_{F_{d-r}}
\frac{\rho^* e_T(U_{d-r})}{e_T(N(F_{d-r}))(- \alpha+c_1(L_{d-r}))}.}
}
Here $e_T(p_{ir}/N_d)=\iota^*_{p_{ir}}\phi_{p_{ir}}$. 
Note that $\varphi^*(\phi_{p_{ir}})$
restricted to $F_r\times F_{d-r}$ is the same as $e_T(p_{ir}/N_d)$ which is
a polynomial only in $ \alpha$ and $\lambda$ as given below.
Similarly for $r=d$, we have 

\eqn\dumbII{
 Q_d(\lambda_i+d \alpha)= \alpha^{-1}e_T(p_{id}/N_d)\sum_{F_d}\int_{F_d}\frac{\rho^* e_T(U_d)}{e_T(N(F_d))( \alpha+c_1(L))}
}
and for $r=0$
\eqn\Qdlambda{
Q_d(\lambda_i)=- \alpha^{-1}e_T(p_{i0}/N_d)\sum_{F_d}\int_{F_d}
\frac{\rho^* e_T(U_d)}{e_T(N(F_d))(- \alpha+c_1(L))} }
We can easily compute $e_T(p_{ir}/N_d)$ which is 
$$e_T(p_{ir}/N_d)=\prod_{(j,m)\neq (i,r)}(\lambda_i-\lambda_j+(r-m) \alpha).$$ 
For $r=0$ and $d$ we have 
$$\eqalign{
e_T(p_{i0}/N_d)&=\prod_{(j,m)\neq (i,0)}(\lambda_i-\lambda_j-m \alpha)
\cr
e_T(p_{id}/N_d)&=\prod_{(j,m)\neq (i,d)}(\lambda_i-\lambda_j+(d-m) \alpha)=\prod_{(j,m)\neq (i,0)}(\lambda_i-\lambda_j+m \alpha).}
$$
The last two identities together with \dumbII, \Qdlambda  clearly gives us 
\eqn\dumbIII{
Q_d(\lambda_i+d \alpha)=\overline{Q_d(\lambda_i)}
}
where $\bar\alpha=-\alpha$, $\bar\lambda_i=\lambda_i$.
Finally our asserted quadratic relation:
\eqn\dumb{
\Omega^V(\lambda_i)Q_d(\lambda_i+r \alpha)=\overline{Q_r(\lambda_i)}
   Q_{d-r}(\lambda_i)}
follows from 
\dumbI, \dumbII, \Qdlambda, \dumbIII, and the following
elementary identity:
$$\prod_{j\neq i}(\lambda_i-\lambda_j)\prod_{(j,m)\neq (i,r)}(\lambda_i-\lambda_j+(r-m) \alpha)=$$
$$\prod_{(j, m)\neq (i,0)}(\lambda_i-\lambda_j-m \alpha)\prod_{(j,m)\neq
(i,0)}(\lambda_i-\lambda_j+m \alpha)$$ 
Note that the last identity is just the interesting identity
$$e_T(T_{p_i}\P^n)\cdot e_T(p_{ir}/N_d)=e_T(p_{id}/N_d)\cdot e_T(p_{i0}/N_d).$$

When $V$ is concave, the fiber at $(f,C)\in\bar\cM_{0,0}(d,\P^n)$ 
of the bundle $U_d$
is $H^1(C,f^*V)$, and we need only one change in the above argument. 
The gluing exact sequence in the concave case is 

$$0\rightarrow V|_{p_i}\rightarrow V_d\rightarrow V_r|_{F_r}\oplus
V_{d-r}|_{F_{d-r}}\rightarrow 0$$ instead. Therefore the gluing identity
for equivariant Euler classes becomes 
$$\chi_d^V=\chi^V_r\cdot \chi^V_{d-r}~\iota^*_{p_i}e_T(V).$$ 
Since $\Omega^V=1/e_T(V)$ for concave $V$, the quadratic relation
\dumb~ remains valid in this case.

The case when $V$ is a direct sum of a convex and a concave bundle
is also similar.  $\Box$

\subsec{Linked Euler data}

\definition{Two sequences $P,Q\in\cS$ are said to  be linked
if $\iota^*_{p_{i,0}}(P_d
-Q_d)\in\cR^{-1}$ 
vanish at $\alpha=(\lambda_i-\lambda_j)/d$
for all $j\neq i$, $d>0$.}

Let $V$ be a $T$-equivariant concavex bundle on $\P^n$, and $C\cong\P^1$ be any
 $T$-invariant line in $\P^n$. By Grothendieck's principle, we have the form
$$V|_C=\oplus_{a=1}^P\cO(l_a)\oplus\oplus_{b=1}^{N}\cO(-k_b)$$
for some positive integers $l_a,k_b$. (0 cannot
occurs because $e_T(V)$ is invertible, by definitiion.) Assume that
$\{l_a\}$ and $\{k_b\}$ are independent of $C$. This is the case, for
example, if $V$ is uniform \Okonek.
{\it We call the numbers $(l_1,..,l_P;k_1,..,k_N)$
the splitting type of $V$.} With this notations, we have

\theorem{Let
$Q_d=\varphi_!(\chi^V_d)$ as before. At 
$\alpha=(\lambda_j-\lambda_i)/d$, $i\neq j$, we have
$$\iota^*_{p_{i,0}}(Q_d)=\prod_a\prod_{m=0}^{l_a d}
(l_a\lambda_j-m(\lambda_j-\lambda_i)/d)\times
\prod_b\prod_{m=0}^{k_b d-1}
(-k_a\lambda_j+m(\lambda_j-\lambda_i)/d).$$
In particular
the Euler data $Q$  is linked to
$$P:~~P_d=
\prod_{a}\prod_{m=0}^{l_a d}(l_a\kappa-m\alpha)
\times
\prod_b\prod_{m=1}^{k_bd-1}(-k_b\kappa+m\alpha).
$$
}
\thmlab\SpecialValueTheorem
\proof
Since we shall evaluate the class $\iota^*_{p_{j,0}}(Q_d)$ at
$\alpha=(\lambda_j-\lambda_i)/d$, we introduce the
notation $Q_d(\kappa,\alpha)=Q_d$, and denote
the value of the class above by 
$Q_d(\lambda_j,(\lambda_j-\lambda_i)/d)$.
First consider the case
$V=\cO(l)$ on $\P^n$. 
We consider a smooth point in $(f, C)\in M_d$ with $C=\P^1$ and in coordinates
$$
f: C\ra \P^1\times \P^n,~~~
[w_0, w_1]\mapsto [w_1, w_0]\times [0,\cdots,w_0^{d}, \cdots,
w_1^{d},\cdots, 0]
$$
where in the last term, $w_0^{d}$ is in the $i$-th
position, $ w_1^{d}$ is in the $j$-th position, and all of the other
components are $0$. The image of $(f, C)$ in $N_d$ under $\varphi$ is the smooth point
$$P_{ij}=[0, \cdots, w_0^d, \cdots, w_1^d,\cdots, 0].$$ 
It is easy to see that, if the weight of $S^1$ in the group $G=S^1\times T$
is $\alpha=(\lambda_j-\lambda_i)/d$, then $P_{ij}$ is fixed by the action
of the subgroup of $G$ with $\alpha=(\lambda_j-\lambda_i)/d$. So $(f,C)$
is a smooth point in $M_d$ fixed by the 
subgroup in $G=S^1\times T$
with $\alpha=(\lambda_j-\lambda_i)/d$. The class
 $Q_d(\kappa, \alpha)$ restricted to $P_{ij}$ is just $Q_d(\lambda_j, (\lambda_j-\lambda_i)/d)$.

At the points $(f, C)\in M_d$ and $P_{ij}\in N_d$, the map $\varphi$ is a canonical identification. From definition, $Q_d(\kappa, \alpha)$ restricted to $P_{ij}$ is
the same as $\chi^V_d$ restricted to $\varphi^{-1}(P_{ij})=(f, C)\in M_d$,
which by definition, is the same as $e_T(U_d)$ restricted to the
$T$-fixed point $(\pi_2\circ f,C)$ in $\bar{\cal M}_{0, 0}(d, \P^n)$. 
Note that $(\pi_2\circ f, C)$ is the degree $d$ cover of the $T$-invariant line 
joining $p_i$ and $p_j$ in $\P^n$. Explicitly 
$$\pi_2\circ f:\ [w_0, w_1]
\rightarrow [0,\cdots, w_0^d,\cdots, w_1^d, \cdots 0].$$ 
So the
induced action of $T$ on $C=\P^1$ has weights $\{ \lambda_i/d, \lambda_j/d\}$.

Now let us compute $e_T(U_d)$ restricted to $(\pi_2\circ f, C)$. When restricted $(\pi_2\circ
f,C)\in\bar{\cal M}_{0, 0}(d, \P^n) $ the fiber of $U_d$ is just the section space $H^0(C, (\pi_2\circ f)^*\cO(l))=H^0(C, \cO(ld))$. It has an explicit basis $\{ w_0^{k}w_1^{ld-k}\}$ with $k=0, \cdots, ld$. Since the $T$-weight of $w_0^{k}w_1^{ld-k}$ is $k\lambda_i/d+(ld-k)\lambda_j/d$, 
by multiplying them together we get
$$Q_d(\lambda_j, (\lambda_j-\lambda_i)/d)
=\prod_{m=0}^{ld}(l\lambda_j-m(\lambda_j-\lambda_i)/d).$$

For a general convex vector bundle $V$ on $\P^n$, 
the fiber of the induced bundle
$V_d$ restricted to $\varphi^{-1}(P_{ij})=(f, C)\in M_d$ is $H^0(\P^1,
(\pi_2\circ f)^*V)$. By  Grothendieck's
 principle, $V$ restricted to the line spanning $p_i,p_j$ splits
into direct sum of line bundles $\{ \cO(l_a) \}$. Pulling them back
to $C$ via the degree $d$ map $\pi_2\circ f$, we get the direct sum
of $\{ \cO(l_ad) \}$.
Since $V$ is convex, each $l_a>0$. By applying the same
computation to each summand, we get  
$$Q_d(\lambda_j,
(\lambda_j-\lambda_i)/d)=\prod_{a}\prod_{m=0}^{l_ad}
(l_a\lambda_j-m(\lambda_j-\lambda_i)/d).$$

For a concave bundle $V$, we need only one
minor change in the above argument. 
We leave to
the reader as an exercise to check that for $V=\cO(-k)$, $k>0$,
$$Q_d(\lambda_j,
(\lambda_j-\lambda_i)/d)=\prod_{m=1}^{kd-1}(-k\lambda_j+m(\lambda_j-\lambda_i)/d),$$
by using either the Atiyah-Bott fixed point formula or by writing down an
explicit basis for $H^1(\P^1, \cO(-kd))$. 
So for a arbitrary concave bundle $V$, we have the form
$$Q_d(\lambda_j,
(\lambda_j-\lambda_i)/d)=\prod_b\prod_{m=1}^{k_bd-1}(-k_b\lambda_j+
m(\lambda_j-\lambda_i)/d).$$
Similarly for an arbitrary concavex bundle $V$, we have the form
$$Q_d(\lambda_j,
(\lambda_j-\lambda_i)/d)=
\prod_{a}\prod_{m=0}^{l_ad}(l_a\lambda_j-m(\lambda_j-\lambda_i)/d)
\times
\prod_b\prod_{m=1}^{k_bd-1}(-k_b\lambda_j+m(\lambda_j-\lambda_i)/d).~~~~~\Box$$

\theorem{Suppose $P,Q$ are any linked 
$\Omega$-Euler data.
If 
$$deg_\alpha~\iota^*_{p_{i,0}}(P_d-Q_d)\leq (n+1)d-2$$
for all $i=0,..,n$ and $d=1,2,..$, then $P=Q$.}
\thmlab\Uniqueness
\proof
By definition, $P_0=Q_0=\Omega$. We will show that $P_d=Q_d$,
assuming that
\eqn\IH{P_r=Q_r,~~r=0,..,d-1.}
Since the $\cR$-valued  pairing $pf_d(u\cdot v)$ on $\cR H^*_G(N_d)$
is nondegenerate, 
it suffices to show that
$$L_s:=pf_d(\kappa^s\cdot (P_d-Q_d))$$
is zero for all $s=0,1,2,...$.
By the localization formula for $pf_d$, we get
$$L_s=
\sum_{i=0}^n\sum_{r=0}^d (\lambda_i+r \alpha)^s
{\iota^*_{p_{ir}}(P_d-Q_d)\over
{\prod^n_{k=0}\prod^d_{m=0}}_{(k,m)\neq(i,r)}
(\lambda_i-\lambda_k+(r-m) \alpha) }.$$
Since $P$ is an Euler data, it follows that
$\iota^*_{p_{ir}}(P_d)$, for each $r=1,..,d-1$, is expressible
in terms of $P_1,..,P_{d-1}$. Likewise for $Q$.
Thus by the inductive hypothesis \IH,
the sum over $r$ above
receives contributions only from the $r=0,d$ terms.
Applying the Reciprocity Lemma (i), we further simplify $L_s$ to 
\eqn\dumbI{\eqalign{
L_s&=\sum_{i=0}^n
\left({\lambda_i^s~ A_i( \alpha)\over  \alpha^d}
+{(\lambda_i+d \alpha)^s~ A_i(- \alpha)\over (- \alpha)^d}\right)\cr
A_i( \alpha)&:=
{(-1)^d\over d!\prod_{k\neq i}(\lambda_i-\lambda_k)}
{\iota^*_{p_{i,0}}(P_d-Q_d)
\over\prod_{k\neq i}\prod_{m=1}^d
(\lambda_i-\lambda_k-m \alpha)}.}
}
Since $P,Q$ are linked Euler data, we have
\eqn\dumb{\iota^*_{p_{i,0}}(P_d-Q_d)=0}
at $ \alpha=(\lambda_i-\lambda_k)/d$,
$k\neq i$.
By the inductive hypothesis \IH~ and the Reciprocity Lemma (iii),
\dumb~holds at $ \alpha=(\lambda_i-\lambda_k)/m$ for $m=1,..,d$ as well.
This shows that $A_i\in\cR=\Q(\lambda)[\alpha]$ for all $i$.
By assumption $deg_ \alpha A_i<(n+1)d-1-nd=d-1$.
But since $L_s\in\cR$ ie. polynomial in $\alpha$, {\it for all $s$},
it follows easily that the $A_i$ must be identically zero.  $\Box$

\subsec{The Lagrange map and mirror transformations}

Throughout this subsection, we fix an invertible class $\Omega$ and
shall denote by $\cA=\cA^\Omega$ the set of $\Omega$-Euler data.

\definition{An invertible map $\mu:\cA\ra\cA$ 
is called a mirror transformation
if for any $P\in\cA$,
$\mu(P)$ is linked to $P$. We call $\mu(P)$ a mirror transform of $P$.}

\definition{(Notations) $\cS_0$ denotes
the set of sequences
$B:~~B_d\in \cR^{-1}H_G^*(N_0),~~d=1,2,...$
We define the map $\cI:\cS\ra\cS_0$, $P\mapsto\cI(P)=B$ where
 $B_d=I^*_d(P_d)$.}

Recall that any equivariant
cohomology class $\omega\in\cR^{-1} H^*_G(N_d)$
is determined by its restrictions
$\iota^*_{p_{ir}}(\omega)\in\cR^{-1}$,
$i=0,..,n$,
$r=0,..,d$.
Conversely given any collection $\omega_{ir}\in\cR^{-1}$, 
there exists a unique
class $\omega\in\cR^{-1}H^*_G(N_d)$ such that
$\iota^*_{p_{ir}}(\omega)=
\omega_{ir}$ for all $i,r$. In fact,
$$\omega=\sum_{i=0}^n\sum_{r=0}^d\omega_{ir}
\prod_{(j,m)\neq(i,r)}{\kappa-\lambda_j-m\alpha\over
\lambda_i-\lambda_j-(m-r)\alpha}.$$
In particular given a sequence $B\in\cS_0$, then for each $d$
there is a unique class $P_d\in\cR^{-1}H^*_G(N_d)$ such that
\eqn\ONE{\iota^*_{p_{ir}}(P_d)=
\iota^*_{p_i}(\Omega)^{-1}~\overline{\iota^*_{p_i}(B_r)}~
\iota^*_{p_i}(B_{d-r}),
~~i=0,..,n,~~ r=0,..,d,}
where we have set $B_0:=\Omega$.
This defines a sequence $P\in\cS$, hence a map
$$\cL_\Omega:\cS_0\ra\cS,~~B\mapsto \cL_\Omega(B)=P.$$
{\it We shall call $\cL=\cL_\Omega$ the Lagrange map.}
By \ONE~ at $r=0$, we get
\eqn\dumb{
\iota^*_{p_i}(B_d)=\iota^*_{p_{i,0}}(P_d)=\iota^*_{p_i}I^*_d(P_d),
~~i=0,..,n.}
First, this implies that the
two classes $B_d,I^*_d(P_d)$ on $N_0=\P^n$ coincide for each $d$.
Thus
\eqn\IL{ B=\cI(P)=\cI\circ\cL(B).}
Thus
$\cL:\cS_0\ra\cS$ is a section of the onto map
$\cI:\cS\ra\cS_0$.
Second,
substituting \dumb~ into \ONE, we get
\eqn\TWO{
\iota^*_{p_i}(\Omega)
\iota^*_{p_{ir}}(P_d)=
\overline{\iota^*_{p_{i,0}}(P_r)}~
\iota^*_{p_{i,0}}(P_{d-r}).}
If, furthermore, we have $P_d\in\cR H_G^*(N_d)$
rather than in $\cR^{-1}H_G^*(N_d)$, then eqn. \TWO~
says that $P$ is an Euler data. 

\item{(A)} {\it The image $P=\cL(B)$ of a given
$B\in\cS_0$ under the Lagrange map
 is an Euler data if $P_d\in \cR H_G^*(N_d)$, $d>0$.}
\bs

On the other hand, it is trivial to show that
if $Q\in\cA\subset\cS$ then
\eqn\LI{Q=\cL\circ\cI(Q).}
Now using $\cL$ we can lift any map
$\mu_0:\cS_0\ra\cS_0$ to a map
$$\mu=\cL\circ \mu_0\circ\cI:\cS\ra\cS,$$
which {\it we shall call the Lagrange lift of $\mu_0$.}
Thus from eqns. \IL~and \LI, we have

\item{(B)} {\it Let $\mu_0:\cS_0\ra\cS_0$ be invertible
with inverse $\nu_0$, and let $\mu,\nu$ be their respective
Lagrange lifts. Then $\mu\circ\nu=\nu\circ\mu=id_\cA$ when
restricted to Euler data.}
\bs

We now discuss the relationship between
Euler data and series of  hypergeometric type.
\definition{
Given any $B\in\cS_0$, define
$$
HG[B](t):=
e^{-pt/\alpha}\left( \Omega+ 
\sum_{d>0}
{B_d~ e^{dt}\over 
\prod_{k=0}^n\prod_{m=1}^d(p-\lambda_k-m\alpha)}\right)
$$
where $p\in H^*_G(N_0)$ is the equivariant hyperplane class
of $N_0=\P^n$.}

Note that $HG[B](t)$ is a cohomology valued formal series. 
If $P:~P_d=\prod^{(n+1)d}_{m=0}(l\kappa-m\alpha)$
as in Example 1, it is obvious that
in the limit $\lambda\ra0$, we have
$$
HG[\cI(P)](t)=e^{-Ht/\alpha} \sum_{d\geq0}
{\prod^{(n+1)d}_{m=0}((n+1)H-m\alpha)\over
\prod_{m=1}^d(H-m\alpha)^{n+1}} e^{dt}
$$
where $H\in H^*(\P^n)$ on the right hand side
is the hyperplane class of $\P^n$.
The coefficients of $(-{H\over\alpha})^i$ for $i=1,..,n$,
are exactly solutions to
a hypergeometric differential equation
discussed in the Introduction.

We now consider a construction of mirror transformations.
Let $B\in\cS_0$, and set $B_0:=\Omega$.
Given any power series $g\in e^t\cR[[e^t]]$,
there is a unique $\tilde B\in\cS_0$
such that
$$HG[B](t+g)
=HG[\tilde B](t).$$
In fact, since
$$
HG[B](t+g)
=e^{-pt/\alpha} e^{-pg/\alpha}\sum_{d\geq0} 
{B_d~ e^{dt} e^{dg}\over 
\prod_{k=0}^n\prod_{m=1}^d(p-\lambda_k-m\alpha)},
$$
if we write $e^{d g}=\sum_{s\geq0} g_{d,s}e^{st}$,
$g_{d,s}\in\cR$ and
$e^{-p g/\alpha}=\sum_{s\geq0} g_s' e^{st}$,
$g_s'\in\cR[p/\alpha]$, then it is straightforward to find that
\eqn\gTransform{\eqalign{
\tilde B_d&=B_d'+\sum_{r=0}^{d-1} g_{d-r}' B_r'\prod_{j=0}^n
\prod_{m=r+1}^d(p-\lambda_j-m\alpha)\cr
B_d'&:=B_d+\sum_{r=0}^{d-1} g_{r,d-r} B_r\prod_{j=0}^n
\prod_{m=r+1}^d(p-\lambda_j-m\alpha).}
}
Thus we have an invertible transformation $\mu_0:\cS_0\ra\cS_0$,
$B\mapsto\tilde B$.
Similarly, given any power series $f\in e^t\cR[[e^t]]$
we have an invertible transformation $\mu_0:\cS_0\ra\cS_0$,
$B\mapsto\tilde B$, such that
$$e^{f/\alpha}~HG[B](t)
=HG[\tilde B](t).$$
Again if we write
$e^{f/\alpha}=\sum_{s\geq0}f_s e^{st}$, $f_s\in\cR[\alpha^{-1}]$,
then
\eqn\fTransform{
\tilde B_d=B_d+\sum_{r=0}^{d-1} f_{d-r} B_r\prod_{j=0}^n
\prod_{m=r+1}^d(p-\lambda_j-m\alpha).
}

We now make an important observation about
the transformation $\mu_0$ in each case above.
For each $r=0,..,d-1$ the class 
$\prod_{j=0}^n(p-\lambda_j-d\alpha)$ 
always vanishes
when restricted to the fixed points $p_i\in\P^n$,
at $\alpha=(\lambda_i-\lambda_j)/d$.
It follows immediately from 
\gTransform~and \fTransform~ that $\iota^*_{p_i}(B_d)$,
$\iota^*_{p_i}(\tilde B_d)$ always agree (whenever defined
for all $d$)
at $\alpha=(\lambda_i-\lambda_j)/d$, for $j\neq i$.
To summarize:

\item{(C)} {\it Given $g,f\in e^t\cR[[e^t]]$, let
$\mu_0:\cS_0\ra\cS_0$, $B\mapsto\tilde B$,
be the invertible transformation defined by}
$$e^{f/\alpha}~ HG[B](t+g)
=HG[\tilde B](t).$$
{\it Suppose $B$ is such that all values of
$\iota^*_{p_i}(B_d)$ are well-defined at
$\alpha=(\lambda_i-\lambda_j)/d$, $j\neq i$.
Then these values are preserved under $\mu_0$.}
\bs

Obviously if $P$ is an Euler data and $B=\cI(P)$, then 
the restrictions
 $\iota^*_{p_i}(B_d)=\iota^*_{p_{i,0}}(P_d)\in\cR$ are
polynomial in $\alpha$. Hence they 
are always well-defined at 
$\alpha=(\lambda_i-\lambda_j)/d$, $j\neq i$.

\lemma{
Let $\mu$ be the Lagrange lift of
the above transformation $\mu_0:\cS_0\ra\cS_0$,
$B\mapsto\tilde B$.
Then $\mu$ is a mirror transformation. In particular, if
$P$ is an Euler data, then $\tilde P=\mu(P)$ is
an Euler data with
$$e^{f/\alpha}~ HG[\cI(P)](t+g)
=HG[\cI(\tilde P)](t).$$}
\thmlab\MirrorTransformation
\proof
The second assertion follows from the first assertion
and the fact that 
$\cI\circ \mu=\cI\circ\cL\circ \mu_0\circ\cI=\mu_0\circ\cI$.

It suffices to consider the two cases $f=0$, $g=0$, separately.
In each case we let $P$ be
an Euler data, and denote
$$\tilde P=\mu(P),~~~B=\cI(P),~~~\tilde B=\mu_0(B).$$
Since $\mu:=\cL\circ \mu_0\circ\cI$, we have $\tilde P=\cL(\tilde B)$.
In each case {\it we will show that
$\tilde P$ is an Euler data.}
We claim that this suffices. First, by statement (B) above,
$\mu$ is invertible as a transformation on the
set $\cA$ of Euler data.
Second, by statement (C) above,
the restrictions $\iota^*_{p_{i,0}}(B_d)=\iota^*_{p_{i,0}}(P_d)$
and $\iota^*_{p_{i,0}}(\tilde B_d)=\iota^*_{p_{i,0}}(\tilde P_d)$
agree at $\alpha=(\lambda_i-\lambda_j)/d$, $j\neq i$.
Thus $\tilde P$ is linked to $P$. So, by definition, $\mu$
is a mirror transformation.

We now proceed to checking Eulerity of $\tilde P$.
Since $P=\cL(B)$, \ONE~ holds.
Multiply both sides of eqn. \ONE~ by  
the respective sides of the following identity:
$$\eqalign{
&e^{d\tau}{e^{(\lambda_i+r\alpha)(t-\tau)/\alpha}\over
{\prod_{j=0}^n\prod_{m=0}^d}_{(j,m)\neq(i,r)}
(\lambda_i+r\alpha-\lambda_j-m\alpha)}\cr
&={1\over\prod_{j\neq i}(\lambda_i-\lambda_j)}\times
e^{\lambda_it/\alpha}{ e^{rt}\over
\prod_{j=0}^n\prod_{m=1}^r(\lambda_i-\lambda_j+m\alpha)}\times
e^{-\lambda_i\tau/\alpha} { e^{(d-r)\tau}\over
\prod_{j=0}^n\prod_{m=1}^{d-r}(\lambda_i-\lambda_j-m\alpha)},
}$$
and then sum over $i=0,..,n$, and $r=0,..,d$. The result is 
$$\eqalign{
&e^{d\tau}~pf_d(P_d~ e^{\kappa(t-\tau)/\alpha})\cr
&=\sum_{r=0}^d
pf\left(\Omega^{-1}~
\overline{[ e^{-pt/\alpha} 
{B_r~e^{rt}\over
\prod_{j=0}^n\prod_{m=1}^r(p-\lambda_j-m\alpha)}]} \times
[e^{-p\tau/\alpha} 
{B_{d-r}~e^{(d-r)\tau}\over
\prod_{j=0}^n\prod_{m=1}^{d-r}(p-\lambda_j-m\alpha)}]\right).}
$$
Now summing this over $d=0,1,2,..$, we get:
\eqn\THREE{
\sum_{d\geq0} e^{d\tau}~pf_d(P_d~ e^{\kappa(t-\tau)/\alpha})
=pf\left(\Omega^{-1}~
\overline{HG[B](t)}~ HG[B](\tau)\right).
}
Likewise, of course, for $\tilde P$ and $\tilde B$.

First case:
\eqn\dumbI{
HG[\tilde B](t)=HG[B](t+g(e^t)).}
By \THREE, we have
$$\eqalign{
pf\left(\Omega^{-1}~
\overline{HG[\tilde B](t)}~ HG[\tilde B](\tau)\right)
&=\sum_{d\geq0} e^{d\tau}~pf_d(\tilde P_d~ e^{\kappa(t-\tau)/\alpha})\cr
pf\left(\Omega^{-1}~
\overline{HG[B](t+g(e^t))}~ HG[B](\tau+g(e^\tau))\right)
&=\sum_{d\geq0} e^{d(\tau+g(e^\tau))}~
pf_d\left(P_d~ e^{\kappa(t+\bar g(e^t)-\tau-g(e^\tau))/\alpha}\right).}
$$
By \dumbI, we can equate the two right hand sides above.
Setting $q=e^\tau$, $\zeta=(t-\tau)/\alpha$, we get
\eqn\dumb{\sum_{d\geq0} q^d~ pf_d(\tilde P_d~ e^{\zeta\kappa})
=\sum_{d\geq0} q^d e^{dg(q)}~
pf_d\left(P_d~ e^{\kappa\zeta}~
e^{\kappa(g_+(q e^{\zeta\alpha})-g_+(q))/\alpha}~
e^{-\kappa(g_-(q e^{\zeta\alpha})+g_-(q))/\alpha}\right)
}
where $g=g_+ +g_-$ with $\bar g_\pm=\pm g_\pm$.
Obviously for any $g(q)\in\cR[[q]]$,
$g_+(q e^{\zeta\alpha})-g_+(q)\in \alpha\cdot\cR[[q,\zeta]]$.
Since the involution $\omega\mapsto\bar\omega$ on $\cR$ simply
changes the sign of $\alpha$, the fact that $g_-$ is odd
shows that $g_-(q)\in \alpha\cdot\cR[[q]]$. Likewise for
$g_-(q e^{\zeta\alpha})$.
We know that $P_d\in \cR H_G^*(N_d)$ (since $P$ is Euler),
and that $pf_d$ maps  $\cR H_G^*(N_d)$ to $\cR$. So
the right hand side of \dumb~ now clearly lies in $\cR[[q,\zeta]]$.
So, likewise for the left hand side of \dumb.
It follows that 
\eqn\dumb{
pf_d(\tilde P_d~\kappa^s)\in\cR,~~s=0,1,2,...}
A priori $\tilde P_d\in\cR^{-1} H^*_G(N_d)$ has the form
$$\tilde P_d=a_N\kappa^N+\cdots+a_0,~~a_i\in\cR^{-1},~~
N=(n+1)d+n.$$
Since $pf_d(\kappa^N)=1$, it follows from \dumb~
that $a_N,..,a_0\in\cR$. Hence $\tilde P_d\in \cR H^*_G(N_d)$
(rather than in $\cR^{-1} H^*_G(N_d)$).
By statement (A) above, $\tilde P$ is an Euler data.

Second case:
$$
HG[\tilde B](t)=e^{f/\alpha}~HG[B](t).
$$
Again applying \THREE~ and writing
 $f\in e^t\cR[[e^t]]$ as $f=f_+ +f_-$ with
$\bar f_\pm=\pm f_\pm$, we get
$$\eqalign{
\sum_{d\geq0} q^d~pf_d(\tilde P_d~ e^{\zeta\kappa})
&=e^{-\bar f(e^t)/\alpha}~
e^{f(e^\tau)/\alpha}~
pf\left(\Omega^{-1}~
\overline{HG[B](t)}~ HG[B](\tau)\right)\cr
&=e^{-(f_+(q e^{\zeta\kappa})-f_+(q))/\alpha}~
e^{(f_-(q e^{\zeta\kappa})+f_-(q))/\alpha}~
\sum_{d\geq0} q^d~pf_d(P_d~ e^{\zeta\kappa}).}
$$
The right hand side lie in $\cR[[q,\zeta]]$ as before,
implying that $\tilde P$ is Euler.  $\Box$

\newsec{Applications}

\definition{A concavex bundle $V$ on $\P^n$ is called
a critical bundle if the induced bundle
$U_d\ra\bar\cM_{0,0}(d,\P^n)$ has rank $(n+1)d+n-3=dim~
\bar\cM_{0,0}(d,\P^n)$.
We denote the nonequivariant Euler class by $c_{top}(U_d)$.}

Given a concavex bundle $V=V^+\oplus V^-$, as before we introduce
$$\eqalign{
\Omega^V&:=e_T(V^+)/e_T(V^-)\cr
Q:~~Q_d&=\varphi_!(\chi^V_d),~~~d>0.}
$$
By Theorem \GeneralEulerTheorem~
the sequence $Q$
is an $\Omega^V$-Euler data.
If $V$ is a critical bundle, introduce
$$
K_d:=
\int_{\bar\cM_{0,0}(d,\P^n)}
c_{top}(U_d),~~~~\Phi:=\sum_{d>0}K_d e^{dt}.
$$

\theorem{Let $V$ be a concavex bundle on $\P^n$.
\item{(i)} The restrictions $I_d^*(Q_d)\in H_G^*(\P^n)$ 
has $deg_\alpha I^*_d(Q_d)\leq(n+1)d-2$.
\item{(ii)} If $V$ is critical, then 
in the nonequivariant limit $\lambda\ra0$,
$$\eqalign{
\int_{\P^n}e^{-Ht/\alpha}
{I^*_d(Q_d)\over\prod_{m=1}^d(H-m\alpha)^{n+1}}
&=\alpha^{-3}(2-d~t)K_d\cr
\int_{\P^n}\left(HG[\cI(Q)](t)-e^{-Ht/\alpha}\Omega^V\right)
&=\alpha^{-3}(2\Phi-t\Phi').
}
$$
}
\thmlab\KdTheorem
\proof
The second equality in (ii) follows trivially from the first
equality. 

By eqn. \Qdlambda~ in the proof of Theorem 
\GeneralEulerTheorem, we have
\eqn\dumb{
Q_d(\lambda_i)=\phi_{p_{i,0}}(\lambda_i)
\sum_{F_d}\int_{F_d}{\rho^* e_T(U_d)\over e_T(N(F_d))[\alpha(\alpha-c_1(L))]}
}
where $\phi_{p_{i,0}}(\lambda_i)=
\phi_{p_i}(\lambda_i)\phi_{\P^n}(\lambda_i)$ with
$\phi_{p_i}:=\prod_{j\neq i}(p-\lambda_j)$,
$\phi_{\P^n}:=\prod_{j=0}^n\prod_{m=1}^d(p-\lambda_j-m\alpha)
\in H_T^*(\P^n)$.
From the localization formula, we deduce
$$\eqalign{
\phi_{p_i}(\lambda_i)\sum_{F_d}\int_{F_d}{\rho^* e_T(U_d)\over e_T(N(F_d))[\alpha(\alpha-c_1(L))]}
&=\int_{\bar\cM_{0,1}(d, \P^n)}{\rho^* e_T(U_d)\over\alpha(\alpha-c_1(L))}ev^*\phi_i\cr
&=\int_{\P^n}ev_!\left({\rho^* e_T(U_d)\over\alpha(\alpha-c_1(L))}\right) \phi_i\cr
&= \iota^*_{p_i}ev_!\left({\rho^* e_T(U_d)\over \alpha(\alpha-c_1(L))}\right).
}
$$
Thus \dumb~ can be written as
$$\iota_{p_i}^*I^*_d(Q_d)=\iota_{p_i}^*[\phi_{\P^n}~
ev_!\left({\rho^* e_T(U_d)\over\alpha(\alpha-c_1(L))}\right)],~~~i=0,..,n.$$
It follows that 
\eqn\dumb{I_d^*(Q_d)=\phi_{\P^n}~
ev_!\left({\rho^* e_T(U_d)\over\alpha(\alpha-c_1(L))}\right).}
This shows that
$deg_\alpha I^*_d(Q_d)\leq
deg_\alpha \phi_{\P^n}-2=(n+1)d-2$, proving (i).

Since $Q_d=\varphi_!(\chi^V_d)\in H_G^*(N_d)$, their
nonequivariant limit $\lambda\ra0$ exist.
In this limit \dumb~ gives
$$\eqalign{
A&:=
\int_{\P^n}e^{-Ht/\alpha}~ 
{I_d^*(Q_d)\over\prod_{m=1}^d(H-m\alpha)^{n+1}} \cr
&=\int_{\bar\cM_{0,1}(d,\P^n)} e^{-ev^*HT/\alpha}
{\rho^* c_{top}(U_d)\over\alpha(\alpha-c_1(L))}\cr 
&=\int_{\bar\cM_{0,0}(d,\P^n)} c_{top}(U_d)~ \rho_!
\left({e^{-ev^*Ht/\alpha}
\over \alpha(\alpha-c_1(L))}\right). }
$$

Now $c_{top}(U_d)$ has degree the same
as the dimension of $\bar \cM_{0,0}(d,\P^n)$.
So, the second factor in the last integrand contributes
a scalar factor given by integration
over a generic fiber $E$ (which is a $\P^1$)
of $\rho$.
So we pick out the degree 1 term in
${e^{-ev^*Ht/\alpha}
\over \alpha(\alpha-c_1(L))}$, which is just
${-ev^*Ht\over \alpha^3}+{c_1(L)\over \alpha^3}$.
Restricting to the generic fiber $E$,
say over $(f,C)\in\cM_{0,0}(d,\P^n)$, the
evaluation map $ev$ is equal to $f$, which is
a degree $d$ map $E\cong\P^1\ra\P^n$.
It follows that
$$\int_E ev^*H=d.$$
Moreover, since $c_1(L)$ restricted to $E$ is just the
first Chern class of the tangent bundle to $E$, it follows that
$$\int_E c_1(L)=2.$$
So we have
$$A=(-{d t\over \alpha^3}+{2\over \alpha^3})K_d.~~~~\Box$$

It is easy to work out the complete list of {\it critical}
concavex vector bundles $V$ on $\P^n$ which are direct sums of
line bundles. Such a $V$ is of the form 
$$\eqalign{
V&=V^+\oplus V^-\cr
V^+&=\oplus_{a=1}^P\cO(l_a)\cr
V^-&=\oplus_{b=1}^N\cO(-k_b)}
$$
where $l_1,..,l_P,~k_1,..,k_N$ are positive integers.
By Riemann-Roch, the obstruction bundles $U_d$ that
$V$ induces on $\bar\cM_{0,0}(d,\P^n)$ has 
$rank~U_d=d\left(\sum l_a+\sum k_b\right)+P-N$,
which must be $(n+1)d+n-3$ for all $d$ if $V$ is critical.
Thus we must have
\eqn\dumb{\eqalign{
\sum l_a+\sum k_b&=n+1\cr
P-N&=n-3.}
}

The complete list of critical bundles that are also direct sums of line
bundles is:
\eqn\CompleteList{
\eqalign{
\P^1:~~~&\cO(-1)\oplus\cO(-1)\cr
\P^2:~~~&\cO(-3)\cr
\P^3:~~~&\cO(2)\oplus\cO(-2)\cr
\P^4:~~~&\cO(5),~~\cO(2)\oplus\cO(2)\oplus\cO(-1)\cr
\P^5:~~~&\cO(2)\oplus\cO(4),~~\cO(3)\oplus\cO(3)\cr
\P^6:~~~&\cO(2)\oplus\cO(2)\oplus\cO(3)\cr
\P^7:~~~&\cO(2)\oplus\cO(2)\oplus\cO(2)\oplus\cO(2).}
}
Note that we have excluded the critical bundles in which
the hyperplane bundle $\cO(1)$ occurs
because in the nonequivariant limit
it only reduces a given case of $\P^n$ to $\P^{n-1}$.
For example, even though
the bundle 
$\cO(1)\oplus\cO(-1)\oplus\cO(-1)$ on $\P^2$
is certainly critical, computing the $K_d$ for
the induced bundles is equivalent to doing the same with
$\cO(-1)\oplus\cO(-1)$ on $\P^1$.
It is curious to note that 
the numerical conditions \dumb~ is rather similar 
to the condition for having a projective complete intersection
Calabi-Yau threefold. In fact five of the examples
on $\P^4$ through $\P^7$ above involving only positive bundles
are exactly the cases  in which each
critical bundle cuts out a complete intersection
Calabi-Yau threefold. We also note that, with
the exception of the $\P^1$ case, the three examples
which involve negative bundles in fact correspond to
{\it noncompact Calabi-Yau threefolds}. The total
space of $\cO(-3)\ra\P^2$, the total space of
$\psi^*\cO(-2)\ra X$, where $\psi:X\hookrightarrow\P^3$
is a quadric, and the total space of
$\psi^*\cO(-1)\ra\P^4$ where $\psi:X\hookrightarrow\P^4$
is the intersection of two quadrics, all three are
noncompact Calabi-Yau. These examples arise in the
so-called {\it local mirror symmetry}.
In the next subsection, we shall compute the Euler classes
of the induced bundles for the list above.

\subsec{The first convex example: The Mirror Conjecture}

Throughout this subsection,  we set $l=n+1$,
consider the convex bundle $V=\cO(l)$ on $\P^n$,
and fix $\Omega^V=lp$.
$P,Q$ shall denote the following two linked Euler data
(cf. Theorems \GeneralEulerTheorem, \SpecialValueTheorem.):
$$\eqalign{
&P:~P_d=\prod_{m=0}^{ld}(l\kappa-m\alpha)\cr
&Q:~Q_d=\varphi_!(\chi^V_d).}
$$

Consider the hypergeometric differential equation
$$\left(({d\over dt})^n-l e^t(l{d\over dt}+1)\cdots(l{d\over dt}+n)\right)h(t)=0.$$
We have seen that a basis 
$f_i$, $i=0,..,n-1$, of solutions can be read off
from the hypergeometric series
(cf. Introduction) in the limit $\lambda\ra0$:
$$
HG[\cI(P)](t)=e^{-Ht/\alpha} \sum_{d\geq0}
{\prod^{ld}_{m=0}(lH-m\alpha)\over
\prod_{m=1}^d(H-m\alpha)^{n+1}} e^{dt}
=lH\left(f_0-f_1{H\over\alpha}+f_2{H^2\over\alpha^2}-\cdots\right).
$$
Recall that
\eqn\dumb{T(t):={f_1\over f_0}=t+{g_1\over f_0}}
is the mirror map of Candelas et al, where
$$
f_0:=\sum_{d\geq0}{(ld)!\over (d!)^{n+1}} e^{dt},~~~~
g_1:=\sum_{d\geq1}{(ld)!\over (d!)^l}
\sum_{m=d+1}^{ld}{l\over m}e^{dt}.
$$

\lemma{
In the limit $\lambda\ra0$, we have
$HG[\cI(Q)](T(t))={1\over f_0} HG[(\cI(P)](t)$.
}
\proof 
By expanding in powers of $\alpha^{-1}$ and using the
assumption that $l=n+1$, we get
\eqn\dumbI{\eqalign{
HG[\cI(P)](t)&=
e^{-pt/\alpha}\sum_{d\geq0}{\prod_{m=0}^{ld}(lp-m\alpha)\over
\prod_{k=0}^n\prod_{m=1}^d(p-\lambda_k-m\alpha)} e^{dt}\cr
&=lp\left[f_0-\alpha^{-1}(p~ f_1+
g_2\sum_{k=0}^n\lambda_k)+\cdots\right]}
}
where 
$g_2=\sum_{d\geq1}{(ld)!\over (d!)^l}
\sum_{m=1}^{d}{1\over m}e^{dt}$.

Put $f:=(log~f_0)\alpha 
+{g_2\over f_0}\sum_{k=0}^n\lambda_k\in e^t\cR[[e^t]]$.
By Lemma \MirrorTransformation, 
we have a mirror transformation $\mu$ such that
\eqn\dumbII{
HG[\cI(\tilde P)](t) 
=e^{f/\alpha}~HG[\cI(P)](t)}
where $\tilde P=\mu(P)$.
Substituting \dumbI~ into \dumbII, we get
\eqn\dumbIII{\eqalign{
HG[\cI(\tilde P)](t) 
&=(1+\alpha^{-1}{g_2\over f_0}\sum_{k=0}^n\lambda_k+\cdots)~f_0^{-1}~
lp\left[f_0-\alpha^{-1}(p~ f_1+
g_2\sum_{k=0}^n\lambda_k)+\cdots\right]\cr
&=lp-\alpha^{-1}lp^2{ f_1\over f_0}+
\cdots.}
}

By Lemma \MirrorTransformation~ again, 
we have a mirror transformation
$\nu$ such that
$$HG[\cI(\tilde Q)](t)=HG[\cI(Q)](t+{g_1\over f_0})$$
where $\tilde Q=\nu(Q)$.
Since, by Theorem \KdTheorem~ (i),
$I^*_d(Q_d)=O(\alpha^{(n+1)d-2})$, it is straightforward to find that
\eqn\dumbIV{
HG[\cI(\tilde Q)](t)
=e^{-p(t+{g_1\over f_0})/\alpha}
(lp+\cdots)
=lp-\alpha^{-1}lp^2 (t+ {g_1\over f_0})+
\cdots.
}
From \dumbIII~ and \dumbIV, we  conclude that for $d>0$,
$$
{I^*_d(\tilde P_d-\tilde Q_d)\over 
\prod_{k=0}^n\prod_{m=1}^d(p-\lambda_k-m\alpha)}\equiv0
$$
modulo order $\alpha^{-2}$, and hence
$$deg_\alpha~\iota^*_{p_{i,0}}(\tilde P_d-\tilde Q_d)\leq (n+1)d-2.$$
But $\tilde P=\mu(P)$ is linked to $P$, and
$\tilde Q=\nu(Q)$ is linked to $Q$.
Since $P$ and $Q$ are linked, it follows that $\tilde P$
and $\tilde Q$ are also linked. By Theorem \Uniqueness,
we have $\tilde P=\tilde Q$. In particular, we have
$$HG[\cI(Q)](T(t))=HG[\cI(\tilde Q)](t)=
HG[\cI(\tilde P)](t)=
e^{f/\alpha}~HG[\cI(P)](t).$$
Since both $P_d,Q_d$ lie in $H_G^*(N_d)$, their 
nonequivariant limit exist.
Taking $\lambda\ra 0$ yields our assertion.  $\Box$

{\it Throughout the rest of this subsection,
we set $l=n+1=5$ and consider the
 critical bundle $\cO(5)\ra\P^4$. 
We assume that we have
taken the nonequivariant limit $\lambda\ra0$.} 
Recall that
$$\eqalign{
F(T)&:={5T^3\over6}+\sum_{d>0}K_d e^{dT}\cr
\cF(T)&:={5\over2}({f_1\over f_0}{f_2\over f_0}-{f_3\over f_0}).}
$$

\theorem{(The Mirror Conjecture) $F=\cF$.}
\thmlab\MirrorConjecture
\proof
Since
\eqn\dumbI{
HG[\cI(P)](t)
=5H\left(f_0- f_1{H\over\alpha}+f_2{H^2\over\alpha^2}-
f_3{H^3\over\alpha^3}\right),}
we will prove that (cf. \CDGP)
\eqn\dumb{
HG[\cI(Q)](T)=5 H \left(1-T {H\over \alpha}+ {F'\over5}{H^2\over\alpha^2} 
-{TF'-2F\over5} {H^3\over\alpha^3}\right).}
Eqns. \dumbI, \dumb~ and the preceding lemma imply $F=\cF$.
Denote the right hand side of \dumb~ by $R$.
Then
$$e^{HT/\alpha}R
=5 H \left(1+ {\Phi'\over5}{H^2\over\alpha^2} 
+{2\Phi\over5} {H^3\over\alpha^3}\right)
=5H+O(e^T)
$$
where $\Phi:=F-{5T^3\over6}=O(e^T)$. Similarly 
$e^{HT/\alpha} HG[\cI(Q)](T)
=5H+O(e^T)$, which also has no {\it polynomial dependence}
on $T$. So \dumb~ is equivalent to
$$
\int_{\P^4}e^{-HT/\alpha}\left(e^{HT/\alpha}HG[\cI(Q)]\right)
=\int_{\P^4}e^{-HT/\alpha}\left(e^{HT/\alpha}R\right).
$$
By Theorem \KdTheorem (ii), this left hand side is
$$
\int_{\P^4}HG[\cI(Q)]
=\alpha^{-3}(2\Phi-T\Phi')+\int_{\P^4} e^{-HT/\alpha}5H
=\alpha^{-3}(2F-TF'),
$$
which coincides with $\int_{\P^4}R$.  $\Box$

It is straightforward to generalize Theorem \MirrorConjecture~
to all other critical 
{\it convex} bundles $V$ in the list \CompleteList. 
In each case, $\Omega^V$ becomes $\prod_a l_a p$
and the Euler data $P$
to be linked with $Q:~~Q_d:=\varphi_!(\chi_d^V)$
is given by $P_d=\prod_a \prod_{m=0}^{l_a d}(l_a\kappa-m\alpha)$.
In the nonequivariant limit, the hypergeometric series
$HG[\cI(P)](t)$ will produce
some hypergeometric functions $f_0,..,f_3$ defining
the function $\cF$.
The generating function $F$ for the $K_d$ is modified by
simply replacing the term ${5T^3\over 6}$
by ${T^3\over 6}\int_X\psi^*H^3$, where $\psi:X\ra\P^n$
is the Calabi-Yau
cut out by $V$.
With these minor modifications in each case, Theorem \MirrorConjecture~
holds.
We leave the details as an exercise for the reader.

\subsec{First concave example: multiple-cover formula}

Let $V$ be the bundle $\cO(-1)\oplus\cO(-1)$ on $\P^1$.
For $d>1$, $V$ induces a rank $2d-2$ bundle
$U_d\ra\bar\cM_{0,0}(d,\P^1)$ whose fiber at $(f,C)$
is the space $H^1(C,f^*V)$, thus $V$ is a critical concave bundle
on $\P^1$. We set $\Omega^V=1/e_T(V)=p^{-2}$.
We shall compute the equivariant classes
$Q_d:=\varphi_!(\chi^V_d)$, and the numbers $K_d$ for this critical bundle. 
Note that by definition $Q_1=1$ and $K_1=1$. 
As a consequence of Theorem \SpecialValueTheorem, 
$$\iota^*_{p_{i,0}}(Q_d)=
\prod_{m=1}^{d-1}(\lambda_i-m(\lambda_i-\lambda_j)/d)^2$$
at $\alpha=(\lambda_i-\lambda_j)/d$, $j\neq i$.
Thus $Q$ is linked to
$$P:~~P_d:=\prod_{m=1}^{d-1}(\kappa-m\alpha)^2,$$
which is a $p^{-2}$-Euler
data. Obviously $deg_\alpha I^*_d(P_d)= 2d-2$.
It follows from Theorem \KdTheorem~(i) that
$deg_\alpha I^*_d(P_d-Q_d)\leq 2d-2$, implying $Q=P$
by Theorem \Uniqueness.

\corollary{$K_d=d^{-3}$.}
\proof
By Theorem \KdTheorem (ii)
in the limit $\lambda\ra0$, we have
$$
\int_{\P^1}
e^{-Ht/\alpha}
{I^*_d(Q_d)\over\prod_{m=1}^d(H-m\alpha)^2}
=\alpha^{-3}(2-d~t)K_d.
$$
Since $Q=P$, we have
$I^*_d(Q_d)=I^*_d(P_d)=\prod^{d-1}_{m=1}(H-m\alpha)^2$, giving 
$$
\int_{\P^1} e^{-Ht/\alpha}
{I^*_d(Q_d)\over
\prod_{m=1}^d(H-m\alpha)^2}
=\alpha^{-3}d^{-3}(2-d~t).~~~~~\Box
$$

\subsec{Second concave example: $K_{\P^2}$}

Let $V$ be the canonical bundle $\cO(-3)\ra\P^2$.
For $d>0$, this bundle induces a rank $3d-1$ bundle
$U_d\ra\bar\cM_{0,0}(d,\P^2)$. Thus
$V$ is a critical concave bundle.
We set $\Omega^V=1/e_T(V)=(-3p)^{-1}$.
We shall compute the equivariant classes
$Q_d:=\varphi_!(\chi^V_d)$, and the numbers $K_d$ for this critical bundle. 
As a consequence of Theorem \SpecialValueTheorem, we have
$$\iota^*_{p_{i,0}}(Q_d)=
\prod_{m=1}^{3d-1}(-3\lambda_i+m(\lambda_i-\lambda_j)/d).$$
at $\alpha=(\lambda_i-\lambda_j)/d$, $j\neq i$.
Thus $Q$ is linked to
$$P:~~P_d:=\prod_{m=1}^{3d-1}(-3\kappa+m\alpha).$$
which is a $(-3p)^{-1}$-Euler
data.

\corollary{$HG[\cI(Q)](t+g)=HG[\cI(P)](t)$ where
$g:=\sum_{d>0}{(-1)^d\over d}{(3d)!\over d!^3} e^{dt}$.}
\proof 
By expanding in powers of $\alpha^{-1}$, we get
\eqn\dumbI{\eqalign{
HG[\cI(P)](t)&=
e^{-pt/\alpha}\left((-3p)^{-1}+\sum_{d>0}
{\prod_{m=1}^{3d-1}(-3p+m\alpha)\over
\prod_{k=0}^2\prod_{m=1}^d(p-\lambda_k-m\alpha)} e^{dt}\right)\cr
&=(-3p)^{-1}+\alpha^{-1}{t+g\over3}+\cdots.}
}
As before, it is now straightforward to show that
$$HG[\cI(Q)](t+g)\equiv HG[\cI(P)](t)$$
modulo order $\alpha^{-2}$. Once again
by Theorem \Uniqueness, the two sides
are equal to all orders. $\Box$

Using Theorem \KdTheorem (ii) and the preceding corollary,
we obtain the $K_d$, for $d=1,..,10$:
 $$        
  \vbox{\offinterlineskip        
  \hrule        
  \halign{ &\vrule# & \strut\quad\hfil#\quad\cr     
  \noalign{\hrule}        
  height1pt&\omit&   &\omit&\cr        
  &$d$&& $K_d$ &\cr        
  \noalign{\hrule}        
  &1&& $3$ &\cr        
  &2&& $-{45 \over8 }$ &\cr        
  &3&& ${244 \over9 }$ &\cr        
  &4&& $-{12333 \over64 }$ &\cr        
  &5&& ${211878 \over125 }$ &\cr        
  &6&& $-{102365 \over6 }$ &\cr        
  &7&& ${64639725 \over343 }$ &\cr        
  &8&& $-{1140830253 \over512 }$ &\cr        
  &9&& ${6742982701 \over243 }$ &\cr        
 &10&& $-{36001193817 \over100 }$ &\cr        
  height1pt&\omit&   &\omit&\cr        
  \noalign{\hrule}}        
  \hrule}        
$$

\subsec{A concavex bundle on $\P^3$.}

Let $V=\cO(2)\oplus\cO(-2)$ on $\P^3$. This is a direct sum
of a convex and a concave bundle.
The induced bundle $U_d\ra\bar\cM_{0,0}(d,\P^3)$,
with fiber at $(f,C)$ being $H^0(C,f^*\cO(2))\oplus H^1(C,f^*\cO(-2))$,
has rank $4d$. We set $\Omega^V=e_T(\cO(2))/e_T(\cO(-2))=-1$.
We shall compute the equivariant classes
$Q_d:=\varphi_!(\chi^V_d)$, 
and the numbers $K_d$ for this critical bundle. 
As a consequence of Theorem \SpecialValueTheorem, we have
$$\iota^*_{p_{i,0}}(Q_d)=
\prod_{m=0}^{2d}(2\lambda_i-m(\lambda_i-\lambda_j)/d)\times
\prod_{m=1}^{2d-1}(-2\lambda_i+m(\lambda_i-\lambda_j)/d).$$
at $\alpha=(\lambda_i-\lambda_j)/d$, $j\neq i$.
Thus $Q$ is linked to
$$P:~~P_d:=\prod_{m=0}^{2d}(2\kappa-m\alpha)\times
\prod_{m=1}^{2d-1}(-2\kappa+m\alpha).$$
which is a $-1$-Euler
data.

\corollary{$HG[\cI(Q)](t+g)=HG[\cI(P)](t)$ where
$g:=\sum_{d>0}{1\over d}{(2d)!^2\over d!^4} e^{dt}$.}
\proof 
By expanding in powers of $\alpha^{-1}$, we get
\eqn\dumbI{\eqalign{
HG[\cI(P)](t)&=
e^{-pt/\alpha}\left(-1+\sum_{d>0}
{\prod_{m=0}^{2d}(2p-m\alpha)\times
\prod_{m=1}^{2d-1}(-2p+m\alpha)\over
\prod_{k=0}^3
\prod_{m=1}^d(p-\lambda_k-m\alpha)} e^{dt}\right)\cr
&=-1+\alpha^{-1}p(t+g)+\cdots,}
}
which, as in the previous examples, agrees with
$HG[\cI(Q)](t+g)$
up to order $\alpha^{-2}$. Hence
Theorem \Uniqueness~ yields our assertion.
$\Box$

Using Theorem \KdTheorem (ii) and the preceding corollary,
we obtain the $K_d$, for $d=1,..,10$:
 $$        
  \vbox{\offinterlineskip        
  \hrule        
  \halign{ &\vrule# & \strut\quad\hfil#\quad\cr     
  \noalign{\hrule}        
  height1pt&\omit&   &\omit&\cr        
  &$d$&& $K_d$ &\cr        
  \noalign{\hrule}        
  &1&& $-4$ &\cr        
  &2&& $-{9 \over2 }$ &\cr        
  &3&& $-{328 \over27 }$ &\cr        
  &4&& $-{777 \over16 }$ &\cr        
  &5&& $-{30004\over125 }$ &\cr        
  &6&& $-{4073\over3 }$ &\cr        
  &7&& $-{2890808 \over343 }$ &\cr        
  &8&& $-{7168777 \over128 }$ &\cr        
  &9&& $-{285797488 \over729 }$ &\cr        
 &10&& $-{714787509 \over250 }$ &\cr        
  height1pt&\omit&   &\omit&\cr        
  \noalign{\hrule}}        
  \hrule}        
$$

\subsec{A concavex bundle on $\P^4$.}

Consider now the critical bundle $V=\cO(2)\oplus\cO(2)\oplus\cO(-1)$ 
on $\P^4$.
The induced bundle $U_d\ra\bar\cM_{0,0}(d,\P^3)$,
has rank $5d+1$. We set $\Omega^V=e_T(\cO(2))^2/e_T(\cO(-1))=-4p$.
We shall compute the equivariant classes
$Q_d:=\varphi_!(\chi^V_d)$, 
and the numbers $K_d$ for this critical bundle. 
As a consequence of Theorem \SpecialValueTheorem, we have
$$\iota^*_{p_{i,0}}(Q_d)=
\prod_{m=0}^{2d}(2\lambda_i-m(\lambda_i-\lambda_j)/d)^2\times
\prod_{m=1}^{d-1}(-\lambda_i+m(\lambda_i-\lambda_j)/d).$$
at $\alpha=(\lambda_i-\lambda_j)/d$, $j\neq i$.
Thus $Q$ is linked to
$$P:~~P_d:=\prod_{m=0}^{2d}(2\kappa-m\alpha)^2\times
\prod_{m=1}^{d-1}(-\kappa+m\alpha).$$
which is a $-4p$-Euler
data.

\corollary{$HG[\cI(Q)](t+g)=HG[\cI(P)](t)$ where
$g:=\sum_{d>0}{(-1)^d\over d}{(2d)!^2\over d!^4} e^{dt}$.}
\proof 
By expanding in powers of $\alpha^{-1}$, we get
\eqn\dumbII{\eqalign{
HG[\cI(P)](t)&=
e^{-pt/\alpha}\left(-4p+\sum_{d>0}
{\prod_{m=0}^{2d}(2p-m\alpha)^2\times
\prod_{m=1}^{d-1}(-p+m\alpha)\over
\prod_{k=0}^4
\prod_{m=1}^d(p-\lambda_k-m\alpha)} e^{dt}\right)\cr
&=-4p+\alpha^{-1}4p^2(t+g)+\cdots,}
}
which, as in the previous examples, agrees with
$HG[\cI(Q)](t+g)$
up to order $\alpha^{-2}$. Hence
Theorem \Uniqueness ~yields our assertion.
$\Box$

We can work out the $K_d$ here as we did before. 
The $K_d$ here can be obtained by taking $K_d$ from the
preceding example on $\P^3$, and multiply it by $4(-1)^d$.
This is so because in
the nonequivariant limit, the hypergeometric
series $HG[\cI(P)](t)$ (cf. \dumbI~ and \dumbII)
in this example on $\P^4$
and the preceding example on $\P^3$
are related by first a multiplication of
$4p$ followed by a change of variable $e^{dt}\mapsto (-1)^d e^{dt}$.

\subsec{General concavex bundles}

In fact the examples above are representative of the most general
concavex bundle. Let $V=V^+\oplus V^-$ be a concavex bundle on $\P^n$,
and let $Q,P$ be as
defined in Theorem \SpecialValueTheorem, and assume that $V$ has splitting type $(l_1,..,l_P;k_1,..,k_N)$. Note that
$\sum l_a+\sum k_b$ is the value of the
class $c_1(V^+)-c_1(V^-)$ on a $T$-invariant $\P^1$ in $\P^n$.

\theorem{ 
If $d\left(\sum l_a+\sum k_b\right)-N\leq (n+1)d-2$ for all $d>0$,
then $Q=P$.
If $d\left(\sum l_a+\sum k_b\right)-N\leq (n+1)d$ for all $d>0$,
then there exists a mirror transformation $\mu$,
depending only on the $l_a,k_b$,
 such that $Q=\mu(P)$.}
\proof
By definition of $P$ in Theorem \SpecialValueTheorem,
$$deg_\alpha I^*_d(P_d)=
d\left(\sum l_a+\sum k_b\right)- N.$$
Consider the first case, where 
this is bounded above by $(n+1)d-2$ for all $d$.
Then by Theorem \KdTheorem, 
$$deg_\alpha (I^*_d(P_d)- I^*_d(Q_d))\leq (n+1)d-2,$$
implying $Q=P$ by Theorem \Uniqueness.

Consider now the second case. Obviously our assumption
implies that 
$\sum l_a+\sum k_b\leq (n+1)$. It is trivial to show that
the only possibilities not covered by the first case are:
(1) $N=0$ and $\sum l_a=n+1$;
(2) $N=1$ and $\sum l_a+k_1=n+1$.
(3) $N=0$ and $\sum l_a=n$;
In each of these cases, a mirror transformation can
be constructed by immitating the previous examples in
a straightforward  way.
Case (1) immitates the example $\cO(5)\ra\P^4$,
while cases (2), (3) immitate the example $\cO(-3)\ra\P^2$.
It is obvious that in each case, the mirror transformation
depends only on the data $l_a,k_b$.
$\Box$

\corollary{Under the same hypotheses as in the preceding theorem,
the Euler data $Q:~~Q_d=\varphi_!(\chi_d^V)$
 depends only on the splitting type,
ie. the numbers $l_a,k_b$, of the concavex bundle $V$ on $\P^n$.}

Note that not every concavex bundle on $\P^n$ is a direct sum
of line bundles. 
For example the tangent bundle is convex, but is not a direct sum
of line bundles.

\subsec{Equivariant total Chern class}

For simplicity, we restrict to convex bundles. Let $V$ be
a rank $r$ convex bundle on $\P^n$, and let
$$c_{tot}(V)=x^r+x^{r-1} c_1(V)+\cdots+c_r(V)$$
be the equivariant Chern polynomial of $V$.
Similarly we denote by $c_{tot}(U_d)$ the equivariant Chern polynomial
for $U_d$. As explained in Example 9,
we can extend the notion of Euler data $Q$ 
allowing $Q_d$ to depend on $x$ polynomially,
simply by replacing the ground field $\Q$ by $\Q[x]$.
Then a similar argument as in
Theorem \GeneralEulerTheorem~ shows that
the sequence $Q_d:=\varphi_!(\pi^* c_{tot}(U_d))$
is also an Euler data in the generalized sense.
Moreover, the analogue of Theorem \SpecialValueTheorem~ holds,
ie. at $\alpha=(\lambda_j-\lambda_i)/d$, we have the form
$$\iota^*_{p_{i,0}}(Q_d)=
\prod_a\prod_{m=0}^{l_ad}(x+l_a\lambda_j-m(\lambda_j-\lambda_i)/d).$$
Hence $Q$ is linked to the Euler data
$$P:~~P_d=
\prod_a\prod_{m=0}^{l_ad}(x+l_a\kappa-m\alpha).$$
Again, under a suitable bound on $c_1(V)$, one can easily
relate $P,Q$ by a generalized (depending on $x$) mirror transformation.

For example, by taking the bundle $\cO(4)$ on $\P^3$,
and applying the above result, we can
compute the nonequivariant limits of
 all $\int_{\bar\cM_{0,0}(d,\P^3)}c_{4d}(U_d)$. They are expected
to count rational curves in a $\P^1$-family of K3 hypersurfaces in
$\P^3$. Similarly we can take $\cO(3)$ on $\P^2$ and compute
$\int_{\bar\cM_{0,0}(d,\P^2)}c_{3d-1}(U_d)$ which should count the number
of rarional curves in a $\P^2$-family of elliptic curves in $\P^2$. Details will be reported in full in our forth-coming paper \LLYII.

\subsec{Concluding remarks}

The most important result we establised
 in this paper is the Mirror Principle.
For simplicity we have restricted our examples in this paper,
to studying only Euler classes and total Chern classes. 
As mentioned in the Introduction,
the Mirror Principle works well for any multiplicative
equivariant characteristic classes. We shall study
in details more examples
of the total Chern class in our forthcoming paper \LLYII.
Generalization to manifolds with torus action will also
be dealt with in details there.

Finally, we make a tantalizing observation which might be of both
physical and mathematical significance.
As we have seen, the set of linked Euler data has
an infinite dimensional transformation group -- the mirror group.
For suitable concavex bundle $V\ra\P^n$,
two special linked Euler data  (cf. Theorem \SpecialValueTheorem)
$Q:~Q_d=\varphi_!(\chi^V_d)$
arising from the nonlinear sigma model (the stable map moduli),
and $P$ the corresponding Euler data  of hypergeometric type,
are related by
a mirror transformation.
Since the mirror group is so big,
there are many other Euler data which are linked to $P$
and can be obtained
simply by acting on $P$ by the mirror group.
From the physical point of view, $P$ arises from type IIB
string theory while $Q$ arises from type IIA string theory,
and mirror symmetry is a duality between the two.
This relationship manifests itself on the linear sigma model
as a duality transformation.
This suggests that some other Euler data
linked to $P$ may
arise from some other string theories which are dual
to type IIA and IIB, via more general mirror transformations. From the point of
view of moduli theory, $P$ is associated
to the linear sigma model compatification $N_d$
of the moduli space $M_d^0$ we discussed in Example 10.
Whereas $Q$ is associated to the nonlinear sigma model $M_d$, which is the stable map compactification of $M_d^0$. This suggests that some other Euler
data linked to $P$ may correspond to other compactifications of $M_d^0$.
If true, we will have an association between
string theories, linked Euler
data, and compactifications
of moduli space of maps, all in the same picture, whereby
there is a duality in each kind which one sees in the linear sigma model.
It would be interesting to understand this duality more
precisely.

\footatend\vfill\supereject\immediate\closeout\rfile\writestoppt
\baselineskip=14pt\centerline{{\bf References}}\bigskip{\frenchspacing%
\parindent=20pt\escapechar=` \input refs.tmp\vfill\eject}\nonfrenchspacing

\end